\documentclass[a4paper,12pt]{article}
\usepackage{times, url}
\usepackage{multirow}

\usepackage{graphicx,amsmath,amssymb,amsthm,amsfonts,bm,float,cite}%,refcheck
\newtheorem{theorem}{Theorem}[section]
\newtheorem{corollary}{Corollary}[section]
\newtheorem{lemma}{Lemma}[section]
\newtheorem{definition}{Definition}[section]
\newtheorem{remark}{Remark}[section]
\usepackage[none]{hyphenat}

\usepackage[center]{caption}
\usepackage[bottom,symbol]{footmisc}
\DeclareFontEncoding{TS1}{}{}
\DeclareFontSubstitution{TS1}{cmr}{m}{n}
\DeclareTextSymbol{\tcrp}{TS1}{'251}
\DeclareTextSymbolDefault{\tcrp}{TS1}

\begin{document}
\setcounter{page}{1}

%\thispagestyle{empty} 
%\noindent \textbf{Notes on Number Theory and Discrete
%Mathematics \newline
%Print ISSN 1310--5132, Online ISSN 2367--8275 \newline
%XXXX, Volume XX, Number X, XXX--XXX \newline
%DOI: 10.7546/nntdm.XXXX}  

\vspace{11mm}

\begin{center}
{\LARGE \bf  The asymptotic estimation for two classes of generalized Fibonacci sub-sequences  } % Sentence case is when the first letter of the title gets capitalized, along with any proper nouns. If your title is 2+ lines use \\[4mm] Please avoid formulas in the title.
\vspace{8mm}

{\Large \bf Yongkang Wan$^1$, Zhonghao Liang$^1$ and Qunying Liao$^{1}$}
\vspace{3mm}
% All authors have their affiliations separately given affiliation, even if they share an affiliation.

$^1$ School of Mathematical Sciences, Sichuan Normal University \\ 
Chengdu, China \\
e-mail: \url{2475636261@qq.com; liangzhongh0807@163.com; qunyingliao@sicnu.edu.cn}
\vspace{2mm}

%$^2$ School of Mathematical Sciences, Sichuan Normal University \\ Chengdu, China \\
%e-mail: \url{liangzhongh0807@163.com}
\vspace{2mm}

%$^3$ School of Mathematical Sciences, Sichuan Normal University \\ 
%Chengdu, China \\
%e-mails: \url{qunyingliao@sicnu.edu.cn}%,  \url{third.author@example.com} % In case you want to indicate 2 (and no more than 2) emails
\end{center}
%\vspace{9mm}  \blfootnote{\licen}

\noindent

%\noindent {\bf Received:}  DD Month XXXX \hfill {\bf Revised:} DD Month XXXX \\
%{\bf Accepted:} DD Month XXXX  \hfill {\bf Online First:} DD Month XXXX \\[4mm] % Please leave these to the Journal's Technical Editor
{\bf Abstract:} Since the $\mathrm{Fibonacci}$ sequence has good properties, it's important in theory and  applications, such as in combinatorics, cryptography, and so on. In this paper, for the generalized Fibonacci sequence $\left\{W_n\left(a,b,p,q\right)\right\}$, by using elementary methods and techniques, we respectively give the asymptotic estimation values of $\left(\sum\limits_{k=n}^{\infty}\frac{1}{W_{mk+l}^d}\right)^{-1}$ and $\left(\sum\limits_{k=n}^{\infty}\frac{\left(-1\right)^k}{W_{mk+l}^d}\right)^{-1}$, which generalize  the asymptotic estimation results of Yuan et al. \cite{A14} in 2025. % Capitalize only first letters of keywords, separate them with commas and finish with  stop.
%{\bf 2020 Mathematics Subject Classification:} XXYZZ, AABCC. % One or more 2020 MSC codes here. Please do note leave this field blank.
\vspace{5mm}

\section{Introduction} \label{sec:Intr}

At the beginning of the 12th century, Fibonacci \cite{A1} proposed the famous Fibonacci sequence $\left\{F_n\right\}$, which is derived from the linear recurrence relation
$$F_0=0,F_1=1,F_n=F_{n-1}+F_{n-2}\left(n\geq 2\right),$$
and the corresponding Binet formula is defined as
$$F_n=\frac{\alpha^n-\beta^n}{\alpha-\beta}\left(n\geq 0\right),$$
where $\alpha=\frac{1+\sqrt{5}}{2}$ and $\beta=\frac{1-\sqrt{5}}{2}$. In 1965, Horadam \cite{A2} defined the generalized Fibonacci sequence $\left\{W_n\left(a,b,p,q\right)\right\}\left(a,b,p,q\in \mathbb{Z}\right)$, which is derived from the linear recurrence relation
$$W_0=a,W_1=b,W_n=pW_{n-1}+qW_{n-2}\left(n\geq 2\right),$$
and the corresponding Binet formula is defined as
$$W_n=c_1\alpha^n-c_2\beta^n\left(n\geq 0\right),$$
where $c_1=\frac{b-a\beta}{\alpha-\beta},c_2=\frac{b-a\alpha}{\alpha-\beta},\alpha=\frac{p+\sqrt{p^2+4q}}{2}$ and $\beta=\frac{p-\sqrt{p^2+4q}}{2}$. 

Since the Fibonacci sequence is important in theory and applications, the various properties of (generalized) Fibonacci sequences have attracted much attention[3-15]. Specifically, in recent years, many scholars considered the asymptotic estimation problem of the reciprocal sums of (generalized) Fibonacci sequences and obtained some results[10-15]. For example, for the Fibonacci sequence $\left\{F_n\right\}$, in 2020, Lee \cite {A11} et al. gave the following two asymptotic estimation values,
$$\left(\sum\limits_{k=n}^{\infty}\frac{1}{F_{k}^2}\right)^{-1}\sim F_{n}^2-F_{n-1}^2+\frac{2}{3}\left(-1\right)^n$$
%$$\left(\sum\limits_{k=n}^{\infty}\frac{1}{F_{k}F_{k+2l}}\right)^{-1}\sim F_{n+l-1}F_{n+l}-\left(F_l^2+\left(-1\right)^l\right)\frac{\left(-1\right)^n}{3},$$
%$$\left(\sum\limits_{k=n}^{\infty}\frac{1}{F_{k}F_{k+2l-1}}\right)^{-1}\sim F_{n+l-1}^2-\left(F_{l-1}F_l+\left(-1\right)^l\right)\frac{\left(-1\right)^n}{3},$$
and
$$\left(\sum\limits_{k=n}^{\infty}\frac{1}{F_{3k}^2}\right)^{-1}\sim F_{3n}^2-F_{3n-3}^2+\frac{4}{9}\left(-1\right)^n.$$
In 2021, Lee \cite {A12} et al. gave the following two asymptotic estimation values,
$$\left(\sum\limits_{k=n}^{\infty}\frac{1}{F_{k}}\right)^{-1}\sim F_{n-2}$$
and
$$\left(\sum\limits_{k=n}^{\infty}\frac{1}{F_{mk-l}}\right)^{-1}\sim F_{mn-l}-F_{m\left(n-1\right)-l},$$
where $m,l$ are positive integers  with $l\leq m-1 $.
In 2022, {Marques\cite{A13} et al. gave the following two asymptotic estimation values,
	$$\left(\sum\limits_{k=2l}^{\infty}\frac{1}{F_{mk}^2}\right)^{-1}\sim F_{2lm}^2-F_{\left(2l-1\right)m}
	^2-\left(-1\right)^{2lm}\frac{2\sqrt{5}\left(L_{2m}-2\right)}{25F_{2m}}$$
	and
	$$\left(\sum\limits_{k=2l+1}^{\infty}\frac{1}{F_{mk}^2}\right)^{-1}\sim F_{\left(2l+1\right)m}^2-F_{2lm}
	^2+\left(-1\right)^{2lm}\frac{2\left(L_{2m}+2\right)}{5L_{2m}},$$
	where $\left\{L_n\right\}$ denotes the Lucas sequence.
	In 2022, Hwang\cite{A10} et al. gave the following asymptotic estimation value,
	$$\left(\sum\limits_{k=n}^{\infty}\frac{1}{F_{k}^4}\right)^{-1}\sim F_{n}^4-F_{n-1}
	^4+\frac{2\left(-1\right)^n}{5}F_{2n-1}+\frac{2\sqrt{5}}{75}.$$
	In 2025, for the generalized $\mathrm{Fibonacci}$ sequence $\left\{W_n\left(a,b,p,q\right)\right\}$, Yuan et al.\cite{A14} firstly considered the corresponding asymptotic estimation problem, and obtained the following four asymptotic estimation results, 
	$$\left(\sum\limits_{k=n}^{\infty}\frac{1}{W_{mk}\left(0,1,A,B\right)}\right)^{-1}\sim W_{mn}\left(0,1,A,B\right)-W_{m\left(n-1\right)}\left(0,1,A,B\right)
	,$$
	$$\left(\sum\limits_{k=n}^{\infty}\frac{1}{W_{mk}^2\left(0,1,A,B\right)}\right)^{-1}\sim W_{mn}^2\left(0,1,A,B\right)-W_{m\left(n-1\right)}^2\left(0,1,A,B\right)
	+B^{mn}C_m,$$
	$$\begin{aligned}
		\left(\sum\limits_{k=n}^{\infty}\frac{1}{W_{mk}^3\left(0,1,A,B\right)}\right)^{-1}\sim& W_{mn}^3\left(0,1,A,B\right)-W_{m\left(n-1\right)}
		^3\left(0,1,A,B\right)\\
		&+3B^{mn}Q_m\left(W_{m\left(n+2\right)}\left(0,1,A,B\right)-W_{m\left(n-3\right)}\left(0,1,A,B\right)\right)
	\end{aligned}$$
	and
	$$\begin{aligned}
		\left(\sum\limits_{k=n}^{\infty}\frac{1}{W_{mk}^4\left(0,1,A,B\right)}\right)^{-1}\sim& W_{mn}^4\left(0,1,A,B\right)-W_{m\left(n-1\right)}
		^4\left(0,1,A,B\right)\\
		&+4B^{mn}U_m\left(W_{m\left(n+1\right)}^2\left(0,1,A,B\right)-B^mW_{m\left(n-2\right)}^2\left(0,1,A,B\right)\right)+V_m,
	\end{aligned}$$
	where 
	$$C_m=\frac{2\left(1-B^m\right)}{\left(\alpha-\beta\right)^2}-\frac{2\left(\alpha^{2m}-1\right)^2}{\left(\alpha-\beta\right)^2\left(\alpha^{4m}-B^m\right)},Q_m=\frac{W_m^2\left(0,1,A,B\right)}{\left(1-\left(B\alpha\right)^{5m}\right)\left(1-\left(B\beta\right)^{5m}\right)},$$
	$$U_m=\frac{W_m^2\left(0,1,A,B\right)}{\left(1-B^m\alpha^{6m}\right)\left(1-B^m\beta^{6m}\right)},V_m=\frac{\left(\alpha^{4m}-1\right)^2}{\left(\alpha-\beta\right)^4}\left(\frac{16\left(\alpha^{4m}-1\right)}{\left(\alpha^{6m}-B^m\right)^2}-\frac{10}{\alpha^{8m}-1}\right).$$
	
	Motivated by the  above works, we 
	consider the most general form of the generalized Fibonacci sequence $\left\{W_n\left(a,b,p,q\right)\right\}$ and obtain two asymptotic estimation results which generalize  the asymptotic estimation results of Yuan et al. \cite{A14} in 2025. 
	
	For the convenience,  we denote $W_n$ by $W_n\left(a,b,p,q\right)$.
	
	This paper is organized as follows. In Section 2, we present the definitions of the asymptotic estimation of convergent series and the Big O notation, and then give some necessary lemmas. In Sections 3-4, we obtain the asymptotic estimation values of $\left(\sum\limits_{k=n}^{\infty}\frac{1}{W_{mk+l}^d}\right)^{-1}$ and $\left(\sum\limits_{k=n}^{\infty}\frac{\left(-1\right)^k}{W_{mk+l}^d}\right)^{-1}$, respectively. In Section 5, we conclude the whole paper.

\section{Preliminaries}

In this section, we review some basic notions and knowledge about the asymptotic estimation of convergent series and the Big O notation, and then give some necessary lemmas. Firstly, we provide the definitions of the asymptotic estimation of convergent series and the Big O notation as the following

%present the definitions of the asymptotic estimation of convergent series, the Big O notation, and some necessary lemmas. Firstly, we provide the definition of the asymptotic estimation of convergent series as the following
\begin{definition}[\cite{A10}]\label{definition 2.1}
	For the convergent series $\sum\limits_{k=1}^{\infty} A_k$, if there exists a function $B_n $ such that
	$$\lim\limits_{n\rightarrow \infty}\left(\left(\sum\limits_{k=n}^{\infty} A_k\right)^{-1}-B_n\right)=0,$$
	then $B_n$ is referred to as the asymptotic estimation of $\left(\sum\limits_{k=n}^{\infty} A_k\right)^{-1}$, and denoted by $\left(\sum\limits_{k=n}^{\infty} A_k\right)^{-1}\sim B_n$.
\end{definition}

%Secondly, we review the definition and some properties of the Big O notation.

\begin{definition}[\cite{A15}]\label{definition 2.2}
	Let $f\left(x\right)$ and $g\left(x\right)\in \mathbb{R}[x]$ with
	$\lim\limits_{x\rightarrow x_0}f\left(x\right)=0$ and $\lim\limits_{x\rightarrow x_0}g\left(x\right)=0$, respectively, where $x_0\in \mathbb{R}$ 
	. If $x$ lies in a certain punctured neighborhood of $x_0 $, and there exists a constant $A>0$ such that
	$$\left|\frac{f\left(x\right)}{g\left(x\right)}\right|\leq A,$$
	then  $\frac{f\left(x\right)}{g\left(x\right)}$ is called a bounded quantity as $x \rightarrow x_0$, and denoted by $f\left(x\right)=O\left(g\left(x\right)\right)\left(x\rightarrow x_0\right).$
\end{definition}

By Definition \ref{definition 2.2}, it's easy to prove  the following Lemmas \ref{lemma 2.1}-\ref{lemma 2.2}.

\begin{lemma}\label{lemma 2.1}
	If $x$ is a real number with $| x |<1 $, then we have
	$$\frac{1}{1+x}=1+\sum_{i=1}^{d-1}\left(-1\right)^ix^i+O\left(x^d\right)$$
	and
	$$\frac{1}{\left(1-x\right)^d}=1+\sum_{i=1}^{d-1}\binom{d-1+i}{d-1} x^i+O\left(x^d\right).$$
\end{lemma}

\begin{lemma}\label{lemma 2.2}
	If $x$ is a real number  with $|x|<1 $, then the following statements are true,
	
	$\left(1 \right)$ for any constant $C$ and positive integer $n$, 
	$$O\left(Cx^n\right)=O\left(x^n\right) \ \text{and}\ C\cdot O\left(x^n\right)=O\left(x^n\right);$$
	
	$\left(2\right)$ for any positive integers $n$ and $m$ with $n \leq m $, 
	$O\left(x^n\right)+O\left(x^m\right)=O\left(x^n\right);$
	
	$\left(3\right)$ for any real number $y$ with $|y|>1 $, positive integers $n$ and $m$, 
	$O\left(\frac{x^n}{y^n}\right)\cdot \frac{1}{y^m}=O\left(\frac{x^n}{y^{n+m}}\right).$
\end{lemma}

%By the proof of Theorem 1 in \cite{A16}, we have the following 
The following Lemma \ref{lemma 2.3} is crucial for the proofs of Theorem \ref{theorem 3.1} and Theorem \ref{theorem 4.1}.

%Finally, to prove our main results, we need the following
\begin{lemma}\label{lemma 2.3}
	For any positive integers $m, k, d, p$ and integers $q,l$, where $l\geq 1-m$, we have
	$$\frac{1}{W_{mk+l}^d}=\frac{1}{c_1^d}\left(\frac{1}{\alpha^{d\left(mk+l\right)}}+\sum_{i=1}^{d-1}\binom{d-1+i}{d-1}\left(\frac{c_2}{c_1}\right)^i\frac{\beta^{i\left(mk+l\right)}}{\alpha^{\left(mk+l\right)\left(d+i\right)}}+O\left(\frac{\beta^{d\left(mk+l\right)}}{\alpha^{2d\left(mk+l\right)}}\right)\right).$$
\end{lemma}
{\bf Proof}.\ By the Binet formula of the generalized Fibonacci sequence $\left\{W_n\left(a,b,p,q\right)\right\}$, it's easy to get
$$
\frac{1}{W_{mk+l}^d}=\frac{1}{\left(c_1\alpha^{mk+l}-c_2\beta^{mk+l}\right)^d}=\frac{1}{c_1^d\alpha^{d\left(mk+l\right)}}\left(1-\frac{c_2}{c_1}\left(\frac{\beta}{\alpha}\right)^{mk+l}\right)^{-d}.
$$
And by $p>0 $, it's easy to prove that $\left|\frac{\beta}{\alpha}\right|<1$. Therefore, for any positive integer $k$, we have
$$\lim\limits_{k\rightarrow \infty}\left(\frac{\beta}{\alpha}\right)^{k}=0,$$
thus by the definition, it's easy to know that there exists a sufficiently large positive integer $k$ such that
$$\left|\frac{c_2}{c_1}\left(\frac{\beta}{\alpha}\right)^{mk+l}\right|<1.$$
Furthermore, by Lemmas \ref{lemma 2.1}-\ref{lemma 2.2}, we have
$$
\begin{aligned}
	\frac{1}{W_{mk+l}^d}&=\frac{1}{c_1^d\alpha^{d\left(mk+l\right)}}\left(1+\sum_{i=1}^{d-1}\binom{d-1+i}{d-1}\left(\frac{c_2}{c_1}\left(\frac{\beta}{\alpha}\right)^{mk+l}\right)^i+O\left(\frac{\beta^{d\left(mk+l\right)}}{\alpha^{d\left(mk+l\right)}}\right)\right)\\
	&=\frac{1}{c_1^d}\left(\frac{1}{\alpha^{d\left(mk+l\right)}}+\sum_{i=1}^{d-1}\binom{d-1+i}{d-1}\left(\frac{c_2}{c_1}\right)^i\frac{\beta^{i\left(mk+l\right)}}{\alpha^{\left(mk+l\right)\left(d+i\right)}}+O\left(\frac{\beta^{d\left(mk+l\right)}}{\alpha^{2d\left(mk+l\right)}}\right)\right).
\end{aligned}
$$

And then we complete the proof of Lemma \ref{lemma 2.3}.

$\hfill\Box$

\section{The asymptotic estimation value of $\left(\sum\limits_{k=n}^{\infty}\frac{1}{W_{mk+l}^d}\right)^{-1}$}
In this section, for the generalized Fibonacci sequence $\left\{W_n\left(a,b,p,q\right)\right\}$, we establish the corresponding connection of $\left(\sum\limits_{k=n}^{\infty}\frac{1}{W_{mk+l}^d}\right)^{-1}$ and its asymptotic estimation(Theorem \ref{theorem 3.1}), and then obtain the more precise asymptotic estimation values of $\left(\sum\limits_{k=n}^{\infty}\frac{1}{W_{mk+l}^d}\right)^{-1}$ for $d=1,2,3,4$, respectively(Corollaries 3.1–3.4).
\subsection{The connection of $\left(\sum\limits_{k=n}^{\infty}\frac{1}{W_{mk+l}^d}\right)^{-1}$ and its asymptotic estimation}
In this subsection, the connection of $\left(\sum\limits_{k=n}^{\infty}\frac{1}{W_{mk+l}^d}\right)^{-1}$ and its asymptotic estimation is given as the following 
\begin{theorem}\label{theorem 3.1}
	For any positive integers $m,n,d,p$ and integers $q,l$, where $l\geq 1-m$  and $p^2+2q-2<p\sqrt{p^2+4q}$, we have
	$$ \begin{aligned}
		&\left(\sum_{k=n}^{\infty}\frac{1}{W_{mk+l}^d}\right)^{-1}\\
		=&W_{mn+l}^d-W_{m\left(n-1\right)+l}^d-\sum_{i=1}^{d}\binom{d}{i}c_1^{d-i}\left(-c_2\right)^{i}\left(\left(\alpha^{mn+l}\right)^{d-i}\left(\beta^{mn+l}\right)^i-\left(\alpha^{m\left(n-1\right)+l}\right)^{d-i}\left(\beta^{m\left(n-1\right)+l}\right)^i\right)\\
		&+\frac{c_1^d\alpha^{d\left(mn+l\right)}\left(\alpha^{md}-1\right)}{\alpha^{md}}\sum_{j=1}^{d-1}\left(-1\right)^j\left(\sum_{i=1}^{d-1}\binom{d-1+i}{d-1}\frac{c_2^{i}\beta^{i\left(mn+l\right)}\alpha^{mi}\left(\alpha^{md}-1\right)}{c_1^{i}\alpha^{i\left(mn+l\right)}\left(\alpha^{m\left(d+i\right)}-\beta^{mi}\right)}\right)^j+O\left(\beta^{d\left(mn+l\right)}\right).
	\end{aligned}$$
	%$$
	%\begin{aligned}
	%&\left(\sum_{k=n}^{\infty}\frac{1}{W_{mk+l}^d}\right)^{-1}\\
	%=&W_{mn+l}^d-W_{m\left(n-1\right)+l}^d-\sum_{i=1}^{d}\binom{d}{i}\left(\left(c_1\alpha^{mn+l}\right)^{d-i}\left(-c_2\beta^{mn+l}\right)^i-\left(c_1\alpha^{m\left(n-1\right)+l}\right)^{d-i}\left(-c_2\beta^{m\left(n-1\right)+l}\right)^i\right)\\
	%&+\frac{c_1^d\alpha^{d\left(mn+l\right)}\left(\alpha^{md}-1\right)}{\alpha^{md}}\sum_{j=1}^{d-1}\left(-1\right)^j\left(\sum_{i=1}^{d-1}\binom{d-1+i}{d-1}\frac{c_2^{i}\beta^{i\left(mn+l\right)}\alpha^{mi}\left(\alpha^{md}-1\right)}{c_1^{i}\alpha^{i\left(mn+l\right)}\left(\alpha^{m\left(d+i\right)}-\beta^{mi}\right)}\right)^j+O\left(\beta^{d\left(mn+l\right)}\right).
	%\end{aligned}
	%$$
\end{theorem}
{\bf Proof}.\ By Lemma \ref{lemma 2.3}, we have
$$
\sum_{k=n}^{\infty}\frac{1}{W_{mk+l}^d}=\frac{1}{c_1^d}\left(\sum_{k=n}^{\infty}\frac{1}{\alpha^{d\left(mk+l\right)}}+\sum_{i=1}^{d-1}\binom{d-1+i}{d-1}\left(\frac{c_2}{c_1}\right)^i\sum_{k=n}^{\infty}\frac{\beta^{i\left(mk+l\right)}}{\alpha^{\left(mk+l\right)\left(d+i\right)}}+\sum_{k=n}^{\infty}O\left(\frac{\beta^{d\left(mk+l\right)}}{\alpha^{2d\left(mk+l\right)}}\right)\right).
$$
Note that
$$\begin{aligned}
	\sum_{k=n}^{\infty}\frac{1}{\alpha^{d\left(mk+l\right)}}&
	=\frac{1}{\alpha^{dl}}\sum_{k=n}^{\infty}\left(\frac{1}{\alpha^{dm}}\right)^k\\
	&=\frac{1}{\alpha^{dl}}\left(\frac{1}{\alpha^{dmn}}+\frac{1}{\alpha^{dmn+dm}}+\frac{1}{\alpha^{dmn+2dm}}+\cdots\right)\\
	&=\frac{1}{\alpha^{dmn+dl}}\left(1+\frac{1}{\alpha^{dm}}+\left(\frac{1}{\alpha^{dm}}\right)^2+\cdots\right)\\
	&=\frac{1}{\alpha^{dmn+dl}}\cdot\frac{1}{1-\frac{1}{\alpha^{dm}}}\\
	&=\frac{\alpha^{md}}{\alpha^{d\left(mn+l\right)}\left(\alpha^{md}-1\right)},
\end{aligned}$$
$$\begin{aligned}
	\sum_{k=n}^{\infty}\frac{\beta^{i\left(mk+l\right)}}{\alpha^{\left(mk+l\right)\left(d+i\right)}}&=\frac{\beta^{li}}{\alpha^{l\left(d+i\right)}}\sum_{k=n}^{\infty}\frac{\beta^{mki}}{\alpha^{mk\left(d+i\right)}}\\
	&=\frac{\beta^{li}}{\alpha^{l\left(d+i\right)}}\left(\frac{\beta^{mni}}{\alpha^{mn\left(d+i\right)}}+\frac{\beta^{mni+mi}}{\alpha^{mn\left(d+i\right)+m\left(d+i\right)}}+\frac{\beta^{mni+2mi}}{\alpha^{mn\left(d+i\right)+2m\left(d+i\right)}}+\cdots\right)\\
	&=\frac{\beta^{li}}{\alpha^{l\left(d+i\right)}}\cdot\frac{\beta^{mni}}{\alpha^{mn\left(d+i\right)}}\left(1+\frac{\beta^{mi}}{\alpha^{m\left(d+i\right)}}+\left(\frac{\beta^{mi}}{\alpha^{m\left(d+i\right)}}\right)^2+\cdots\right)\\
	&=\frac{\beta^{i\left(mn+l\right)}}{\alpha^{\left(mn+l\right)\left(d+i\right)}}\cdot\frac{1}{1-\frac{\beta^{mi}}{\alpha^{m\left(d+i\right)}}}\\
	&=\frac{\beta^{i\left(mn+l\right)}\alpha^{m\left(d+i\right)}}{\alpha^{\left(mn+l\right)\left(d+i\right)}\left(\alpha^{m\left(d+i\right)}-\beta^{mi}\right)},
\end{aligned}$$
and by Lemma \ref{lemma 2.2}, it's easy to get
$$\sum_{k=n}^{\infty}O\left(\frac{\beta^{d\left(mk+l\right)}}{\alpha^{2d\left(mk+l\right)}}\right)=O\left(\frac{\beta^{d\left(mn+l\right)}}{\alpha^{2d\left(mn+l\right)}}\right),$$
thus
$$\small \begin{aligned}
	&\left(\sum_{k=n}^{\infty}\frac{1}{W_{mk+l}^d}\right)^{-1}\\
	=&\left(\frac{1}{c_1^d}\left(\frac{\alpha^{md}}{\alpha^{d\left(mn+l\right)}\left(\alpha^{md}-1\right)}+\sum_{i=1}^{d-1}\binom{d-1+i}{d-1}\frac{c_2^i\beta^{i\left(mn+l\right)}\alpha^{m\left(d+i\right)}}{c_1^i\alpha^{\left(mn+l\right)\left(d+i\right)}\left(\alpha^{m\left(d+i\right)}-\beta^{mi}\right)}+O\left(\frac{\beta^{d\left(mn+l\right)}}{\alpha^{2d\left(mn+l\right)}}\right)\right)\right)^{-1}\\
	=&\left(\frac{\alpha^{md}}{c_1^d\alpha^{d\left(mn+l\right)}\left(\alpha^{md}-1\right)}\left(1+\sum_{i=1}^{d-1}\binom{d-1+i}{d-1}\frac{c_2^i\beta^{i\left(mn+l\right)}\alpha^{mi}\left(\alpha^{md}-1\right)}{c_1^i\alpha^{i\left(mn+l\right)}\left(\alpha^{m\left(d+i\right)}-\beta^{mi}\right)}+O\left(\frac{\beta^{d\left(mn+l\right)}}{\alpha^{d\left(mn+l\right)}}\right)\right)\right)^{-1}\\
	=&\frac{c_1^d\alpha^{d\left(mn+l\right)}\left(\alpha^{md}-1\right)}{\alpha^{md}}\left(1+\sum_{i=1}^{d-1}\binom{d-1+i}{d-1}\frac{c_2^i\beta^{i\left(mn+l\right)}\alpha^{mi}\left(\alpha^{md}-1\right)}{c_1^i\alpha^{i\left(mn+l\right)}\left(\alpha^{m\left(d+i\right)}-\beta^{mi}\right)}+O\left(\frac{\beta^{d\left(mn+l\right)}}{\alpha^{d\left(mn+l\right)}}\right)\right)^{-1}.
\end{aligned}$$
And by $p>0$, it's easy to prove that $\left|\frac{\beta}{\alpha}\right|<1$. Thus, for any positive integer $n$, we have
$$\lim_{n\rightarrow \infty}\left(\sum_{i=1}^{d-1}\binom{d-1+i}{d-1}\frac{c_2^i\beta^{il}\alpha^{mi}\left(\alpha^{md}-1\right)}{c_1^i\alpha^{il}\left(\alpha^{m\left(d+i\right)}-\beta^{mi}\right)}\cdot\left(\frac{\beta}{\alpha}\right)^{imn}+O\left(\frac{\beta^{d\left(mn+l\right)}}{\alpha^{d\left(mn+l\right)}}\right)\right)=0,$$
then by the definition, it's easy to know that there exists a sufficiently large positive integer $n$ such that
$$\left|\sum_{i=1}^{d-1}\binom{d-1+i}{d-1}\frac{c_2^i\beta^{i\left(mn+l\right)}\alpha^{mi}\left(\alpha^{md}-1\right)}{c_1^{i}\alpha^{i\left(mn+l\right)}\left(\alpha^{m\left(d+i\right)}-\beta^{mi}\right)}+O\left(\frac{\beta^{d\left(mn+l\right)}}{\alpha^{d\left(mn+l\right)}}\right)\right|<1.$$
Furthermore, by Lemma \ref{lemma 2.1}, we have
$$\small \begin{aligned}
	&\left(\sum_{k=n}^{\infty}\frac{1}{W_{mk+l}^d}\right)^{-1}\\
	=&\frac{c_1^d\alpha^{d\left(mn+l\right)}\left(\alpha^{md}-1\right)}{\alpha^{md}}\left(1+\sum_{j=1}^{d-1}\left(-1\right)^j\left(\sum_{i=1}^{d-1}\binom{d-1+i}{d-1}\frac{c_2^{i}\beta^{i\left(mn+l\right)}\alpha^{mi}\left(\alpha^{md}-1\right)}{c_1^{i}\alpha^{i\left(mn+l\right)}\left(\alpha^{m\left(d+i\right)}-\beta^{mi}\right)}\right)^j+O\left(\frac{\beta^{d\left(mn+l\right)}}{\alpha^{d\left(mn+l\right)}}\right)\right)\\
	=&c_1^d\alpha^{dmn+dl}-c_1^d\alpha^{dm\left(n-1\right)+dl}+\frac{c_1^d\alpha^{d\left(mn+l\right)}\left(\alpha^{md}-1\right)}{\alpha^{md}}\sum_{j=1}^{d-1}\left(-1\right)^j\left(\sum_{i=1}^{d-1}\binom{d-1+i}{d-1}\frac{c_2^{i}\beta^{i\left(mn+l\right)}\alpha^{mi}\left(\alpha^{md}-1\right)}{c_1^{i}\alpha^{i\left(mn+l\right)}\left(\alpha^{m\left(d+i\right)}-\beta^{mi}\right)}\right)^j\\
	&+O\left(\beta^{d\left(mn+l\right)}\right).
\end{aligned}$$
Note that for any positive integer $e$, we have
$$W_{me+l}^d=c_1^d\alpha^{dme+dl}+\sum_{i=1}^{d}\binom{d}{i}c_1^{d-i}\left(-c_2\right)^{i}\left(\alpha^{me+l}\right)^{d-i}\left(\beta^{me+l}\right)^i,$$
thus
$$ \begin{aligned}
	&\left(\sum_{k=n}^{\infty}\frac{1}{W_{mk+l}^d}\right)^{-1}\\
	=&W_{mn+l}^d-W_{m\left(n-1\right)+l}^d-\sum_{i=1}^{d}\binom{d}{i}c_1^{d-i}\left(-c_2\right)^{i}\left(\left(\alpha^{mn+l}\right)^{d-i}\left(\beta^{mn+l}\right)^i-\left(\alpha^{m\left(n-1\right)+l}\right)^{d-i}\left(\beta^{m\left(n-1\right)+l}\right)^i\right)\\
	&+\frac{c_1^d\alpha^{d\left(mn+l\right)}\left(\alpha^{md}-1\right)}{\alpha^{md}}\sum_{j=1}^{d-1}\left(-1\right)^j\left(\sum_{i=1}^{d-1}\binom{d-1+i}{d-1}\frac{c_2^{i}\beta^{i\left(mn+l\right)}\alpha^{mi}\left(\alpha^{md}-1\right)}{c_1^{i}\alpha^{i\left(mn+l\right)}\left(\alpha^{m\left(d+i\right)}-\beta^{mi}\right)}\right)^j+O\left(\beta^{d\left(mn+l\right)}\right).
\end{aligned}$$

This completes the proof of Theorem \ref{theorem 3.1}.

$\hfill\Box$

\subsection{Some Corollaries}
In this subsection, for some special $d$, we can obtain more precise asymptotic estimate values of $\left(\sum\limits_{k=n}^{\infty}\frac{1}{W_{mk+l}^d}\right)^{-1}$ as the following Corollaries \ref{corollary 3.1}-\ref{corollary 3.4}.

\begin{corollary}\label{corollary 3.1}
	For any positive integers $m,n,p$ and integers $q,l$, where $l\geq 1-m$  and $p^2+2q-2<p\sqrt{p^2+4q}$, we have
	$$\left(\sum_{k=n}^{\infty}\frac{1}{W_{mk+l}}\right)^{-1}\sim W_{mn+l}-W_{m\left(n-1\right)+l}.$$
\end{corollary}
{\bf Proof}.\ By Theorem \ref{theorem 3.1}, we have
$$\left(\sum_{k=n}^{\infty}\frac{1}{W_{mk+l}}\right)^{-1}=W_{mn+l}-W_{m\left(n-1\right)+l}+c_2\beta^{m\left(n-1\right)+l}\left(\beta^m-1\right)+O\left(\beta^{mn+l}\right).$$
And by $p^2+2q-2<p\sqrt{p^2+4q}$, it's easy to prove that $|\beta|<1$, hence,
$$
\begin{aligned}
	&\lim_{n\rightarrow \infty}
	\left(\left(\sum_{k=n}^{\infty}\frac{1}{W_{mk+l}}\right)^{-1}-\left(W_{mn+l}-W_{m\left(n-1\right)+l}\right)\right)\\
	=&\lim_{n\rightarrow \infty}
	\left(c_2\beta^{m\left(n-1\right)+l}\left(\beta^m-1\right)+O\left(\beta^{mn+l}\right)\right)\\
	=&0.
\end{aligned}
$$
Therefore, by Definition \ref{definition 2.1}, we get
$$\left(\sum_{k=n}^{\infty}\frac{1}{W_{mk+l}}\right)^{-1}\sim W_{mn+l}-W_{m\left(n-1\right)+l}.$$

This completes the proof of Corollary \ref{corollary 3.1}.

$\hfill\Box$

\begin{remark}
	By taking $a=0,b=1,p=A\in \mathbb{Z}^+,q=\pm 1$ and $l=0$  in Corollary \ref{corollary 3.1}, the corresponding  result is just Theorem 2.1 in \cite{A14}.
\end{remark}

\begin{corollary}\label{corollary 3.2}
	For any positive integers $m,n,p$ and integers $q,l$, where $l\geq 1-m$  and $p^2+2q-2<p\sqrt{p^2+4q}$, we have
	$$\left(\sum_{k=n}^{\infty}\frac{1}{W_{mk+l}^2}\right)^{-1}\\
	\sim W_{mn+l}^2-W_{m\left(n-1\right)+l}^2+2c_1c_2\left(-q\right)^{mn+l}\left(1-\left(-q\right)^{-m}-\frac{\left(\alpha^{2m}-1\right)^2}{\alpha^m\left(\alpha^{3m}-\beta^m\right)}\right).$$
\end{corollary}
{\bf Proof}.\ By Theorem \ref{theorem 3.1}, we have
$$\begin{aligned}
	&\left(\sum_{k=n}^{\infty}\frac{1}{W_{mk+l}^2}\right)^{-1}\\
	=&W_{mn+l}^2-W_{m\left(n-1\right)+l}^2-\sum_{i=1}^{2}\binom{2}{i}c_1^{2-i}\left(-c_2\right)^{i}\left(\left(\alpha^{mn+l}\right)^{2-i}\left(\beta^{mn+l}\right)^i-\left(\alpha^{m\left(n-1\right)+l}\right)^{2-i}\left(\beta^{m\left(n-1\right)+l}\right)^i\right)\\
	&-\frac{c_1^2\alpha^{2\left(mn+l\right)}\left(\alpha^{2m}-1\right)}{\alpha^{2m}}\cdot \frac{2c_2\beta^{mn+l}\alpha^{m}\left(\alpha^{2m}-1\right)}{c_1\alpha^{mn+l}\left(\alpha^{3m}-\beta^{m}\right)}+O\left(\beta^{2\left(mn+l\right)}\right)\\
	=&W_{mn+l}^2-W_{m\left(n-1\right)+l}^2+2c_1c_2\left(-q\right)^{mn+l}\left(1-\left(-q\right)^{-m}-\frac{\left(\alpha^{2m}-1\right)^2}{\alpha^m\left(\alpha^{3m}-\beta^m\right)}\right)-c_2^2\beta^{2\left(m\left(n-1\right)+l\right)}\left(\beta^{2m}-1\right)\\
	&+O\left(\beta^{2\left(mn+l\right)}\right).
\end{aligned}$$
And by $p^2+2q-2<p\sqrt{p^2+4q}$, it's easy to prove that $|\beta|<1$, hence,
$$\small
\begin{aligned}
	&\lim_{n\rightarrow \infty}
	\left(\left(\sum_{k=n}^{\infty}\frac{1}{W_{mk+l}^2}\right)^{-1}-\left(W_{mn+l}^2-W_{m\left(n-1\right)+l}^2+2c_1c_2\left(-q\right)^{mn+l}\left(1-\left(-q\right)^{-m}-\frac{\left(\alpha^{2m}-1\right)^2}{\alpha^m\left(\alpha^{3m}-\beta^m\right)}\right)\right)\right)\\
	&=\lim_{n\rightarrow \infty}
	\left(-c_2^2\beta^{2\left(m\left(n-1\right)+l\right)}\left(\beta^{2m}-1\right)+O\left(\beta^{2\left(mn+l\right)}\right)\right)\\
	&=0.
\end{aligned}
$$
Therefore, by Definition \ref{definition 2.1}, we can get
$$\left(\sum_{k=n}^{\infty}\frac{1}{W_{mk+l}^2}\right)^{-1}\\
\sim W_{mn+l}^2-W_{m\left(n-1\right)+l}^2+2c_1c_2\left(-q\right)^{mn+l}\left(1-\left(-q\right)^{-m}-\frac{\left(\alpha^{2m}-1\right)^2}{\alpha^m\left(\alpha^{3m}-\beta^m\right)}\right).$$

This completes the proof of Corollary \ref{corollary 3.2}.

$\hfill\Box$

\begin{remark}
	By taking $a=0,b=1,p=A\in \mathbb{Z}^+,q=\pm 1$ and $l=0$  in Corollary \ref{corollary 3.2}, the corresponding  result is just Theorem 2.2 in \cite{A14}.
\end{remark}

\begin{corollary}\label{corollary 3.3}
	For any positive integers $m,n,p$ and integers $q,l$, where $l\geq 1-m$  and $p^2+2q-2<p\sqrt{p^2+4q}$, we have
	$$\left(\sum_{k=n}^{\infty}\frac{1}{W_{mk+l}^3}\right)^{-1}\sim W_{mn+l}^3-W_{m\left(n-1\right)+l}^3+C+D,$$
	where
	$$C=3c_1^2c_2\left(-\alpha q\right)^{mn+l}\left(1-\left(-\alpha q\right)^{-m}-\frac{\left(\alpha^{3m}-1\right)^2}{\alpha^{2m}\left(\alpha^{4m}-\beta^m\right)}\right)$$
	and 
	$$D=3c_1c_2^2\left(-\beta q\right)^{mn+l}\left(-1+\left(-\beta q\right)^{-m}-\frac{2\left(\alpha^{3m}-1\right)^2}{\alpha^m\left(\alpha^{5m}-\beta^m\right)}+\frac{3\left(\alpha^{3m}-1\right)^3}{\alpha^m\left(\alpha^{4m}-\beta^m\right)^2}\right).$$
\end{corollary}
{\bf Proof}.\ By Theorem \ref{theorem 3.1}, we have
$$\begin{aligned}
	&\left(\sum_{k=n}^{\infty}\frac{1}{W_{mk+l}^3}\right)^{-1}\\
	=&W_{mn+l}^3-W_{m\left(n-1\right)+l}^3-\sum_{i=1}^{3}\binom{3}{i}c_1^{3-i}\left(-c_2\right)^{i}\left(\left(\alpha^{mn+l}\right)^{3-i}\left(\beta^{mn+l}\right)^i-\left(\alpha^{m\left(n-1\right)+l}\right)^{3-i}\left(\beta^{m\left(n-1\right)+l}\right)^i\right)\\
	&+\frac{c_1^3\alpha^{3\left(mn+l\right)}\left(\alpha^{3m}-1\right)}{\alpha^{3m}}\sum_{j=1}^{2}\left(-1\right)^j\left(\sum_{i=1}^{2}\binom{2+i}{2}\frac{c_2^{i}\beta^{i\left(mn+l\right)}\alpha^{mi}\left(\alpha^{3m}-1\right)}{c_1^{i}\alpha^{i\left(mn+l\right)}\left(\alpha^{m\left(3+i\right)}-\beta^{mi}\right)}\right)^j+O\left(\beta^{3\left(mn+l\right)}\right)\\
	=&W_{mn+l}^3-W_{m\left(n-1\right)+l}^3+C+D+c_2^3\beta^{3\left(m\left(n-1\right)+l\right)}\left(\beta^{3m}-1+\frac{36\alpha^m\beta^{3m}\left(\alpha^{3m}-1\right)^3}{\left(\alpha^{4m}-\beta^{m}\right)\left(\alpha^{5m}-\beta^{2m}\right)}\right)\\
	&+\frac{36c_2^4\beta^{4\left(mn+l\right)}\alpha^m\left(\alpha^{3m}-1\right)^3}{c_1\alpha^{mn+l}\left(\alpha^{5m}-\beta^{2m}\right)^2}+O\left(\beta^{3\left(mn+l\right)}\right),
\end{aligned}$$
where
$$C=3c_1^2c_2\left(-\alpha q\right)^{mn+l}\left(1-\left(-\alpha q\right)^{-m}-\frac{\left(\alpha^{3m}-1\right)^2}{\alpha^m\left(\alpha^{4m}-\beta^m\right)^2}\right),$$
and
$$D=3c_1c_2^2\left(-\beta q\right)^{mn+l}\left(-1+\left(-\beta q\right)^{-m}-\frac{2\left(\alpha^{3m}-1\right)^2}{\alpha^m\left(\alpha^{5m}-\beta^m\right)}+\frac{3\left(\alpha^{3m}-1\right)^3}{\alpha^m\left(\alpha^{4m}-\beta^m\right)^2}\right).$$
By $p>0$, it's easy to prove that $|\alpha|>|\beta|$. And by  $p^2+2q-2<p\sqrt{p^2+4q}$, it's easy to prove that $|\beta|<1$. Then, $|\alpha|>|\beta|>|\beta^4|$. Furthermore,
$$
\footnotesize
\begin{aligned}
	&\lim_{n\rightarrow \infty}
	\left(\left(\sum_{k=n}^{\infty}\frac{1}{W_{mk+l}^3}\right)^{-1}-\left(W_{mn+l}^3-W_{m\left(n-1\right)+l}^3+C+D\right)\right)\\
	%=&\lim_{n\rightarrow \infty}
	%\left(c_2^3\beta^{3\left(mn+l\right)}-c_2^3\beta^{3\left(m\left(n-1\right)+l\right)}+\frac{36c_2^3\beta^{3\left(mn+l\right)}\alpha^m\left(\alpha^{3m}-1\right)^3}{\left(\alpha^{4m}-\beta^{m}\right)\left(\alpha^{5m}-\beta^{2m}\right)}+\frac{36c_2^4\beta^{4\left(mn+l\right)}\alpha^m\left(\alpha^{3m}-1\right)^3}{c_1\alpha^{mn+l}\left(\alpha^{5m}-\beta^{2m}\right)^2}\right.\\
	%&\left. +O\left(\beta^{3\left(mn+l\right)}\right)\right)\\
	=&\lim_{n\rightarrow \infty}
	\left(c_2^3\beta^{3\left(m\left(n-1\right)+l\right)}\left(\beta^{3m}-1+\frac{36\alpha^m\beta^{3m}\left(\alpha^{3m}-1\right)^3}{\left(\alpha^{4m}-\beta^{m}\right)\left(\alpha^{5m}-\beta^{2m}\right)}\right)+\left(\frac{\beta^4}{\alpha}\right)^{mn+l}\cdot \frac{36c_2^4\alpha^m\left(\alpha^{3m}-1\right)^3}{c_1\left(\alpha^{5m}-\beta^{2m}\right)^2}+O\left(\beta^{3\left(mn+l\right)}\right)\right)\\
	=&0.
\end{aligned}
$$
Therefore, by Definition \ref{definition 2.1}, we can get
$$\left(\sum_{k=n}^{\infty}\frac{1}{W_{mk+l}^3}\right)^{-1}\sim W_{mn+l}^3-W_{m\left(n-1\right)+l}^3+C+D,$$
where
$$C=3c_1^2c_2\left(-\alpha q\right)^{mn+l}\left(1-\left(-\alpha q\right)^{-m}-\frac{\left(\alpha^{3m}-1\right)^2}{\alpha^{2m}\left(\alpha^{4m}-\beta^m\right)}\right)$$
and 
$$D=3c_1c_2^2\left(-\beta q\right)^{mn+l}\left(-1+\left(-\beta q\right)^{-m}-\frac{2\left(\alpha^{3m}-1\right)^2}{\alpha^m\left(\alpha^{5m}-\beta^m\right)}+\frac{3\left(\alpha^{3m}-1\right)^3}{\alpha^m\left(\alpha^{4m}-\beta^m\right)^2}\right).$$

This completes the proof of Corollary \ref{corollary 3.3}.

$\hfill\Box$

The following Remark 3.3 shows that the corresponding result is just Theorem 2.3 in \cite{A14} when $a=0,b=1,p=A\in \mathbb{Z}^+,q=\pm 1$ and $l=0$ in Corollary \ref{corollary 3.3}.
\begin{remark}
	By taking $a=0,b=1,p=A\in \mathbb{Z}^+,q=\pm 1$ and $l=0$ in Corollary \ref{corollary 3.3}, we have $c_1=c_2=\frac{1}{\alpha-\beta}$, and then $$C=\frac{3\left(B\alpha\right)^{mn}}{\left(\alpha-\beta\right)^3}\left(1-B^{-m}\alpha^{-m}-\frac{\left(\alpha^{3m}-1\right)^2}{\alpha^m\left(\alpha^{5m}-\alpha^m\beta^m\right)}\right),$$
	and 
	$$D=\frac{3B^{mn}\beta^{mn}}{\left(\alpha-\beta\right)^3}\left(-1+B^{-m}\beta^{-m}-\frac{2\left(\alpha^{3m}-1\right)^2}{\alpha^m\left(\alpha^{5m}-\beta^m\right)}+\frac{3\left(\alpha^{3m}-1\right)^3}{\alpha^m\left(\alpha^{4m}-\beta^m\right)^2}\right).$$
	Note that $\alpha\beta=B$, it's easy to know that $\alpha^{-m}=B^{-m}
	\beta^m$, and so
	$$C=\frac{3\left(B\alpha\right)^{mn}}{\left(\alpha-\beta\right)^3}\left(1-\beta^m-\frac{\left(\alpha^{3m}-1\right)^2}{\alpha^m\left(\alpha^{5m}-B^m\right)}\right).$$
	Now by the proof of Theorem 2.3 in \cite{A14}, we have  
	$$C=3B^{mn}Q_m\left(W_{m\left(n+2\right)}-W_{m\left(n-3\right)}+\frac{\beta^{m\left(n+2\right)}-\beta^{m\left(n-3\right)}}{\alpha-\beta}\right),$$
	where $Q_m=\frac{W_m^2}{\left(1-\left(B\alpha\right)^{5m}\right)\left(1-\left(B\beta\right)^{5m}\right)}$. 
	And by $\left|\beta\right|<1$, it's easy to get
	$$
	\begin{aligned}
		&\lim_{n\rightarrow \infty}
		\left(\left(\sum_{k=n}^{\infty}\frac{1}{W_{mk}^3}\right)^{-1}-\left(W_{mn}^3-W_{m\left(n-1\right)}^3+3B^{mn}Q_m\left(W_{m\left(n+2\right)}-W_{m\left(n-3\right)}\right)\right)\right)\\
		=&\lim_{n\rightarrow \infty}
		\left(3B^{mn}Q_m\cdot \frac{\beta^{m\left(n+2\right)}-\beta^{m\left(n-3\right)}}{\alpha-\beta}+D\right)\\
		=&0.
	\end{aligned}
	$$
	Therefore, by Definition \ref{definition 2.1}, we have
	$$\left(\sum_{k=n}^{\infty}\frac{1}{W_{mk}^3}\right)^{-1}\sim W_{mn}^3-W_{m\left(n-1\right)}^3+3B^{mn}Q_m\left(W_{m\left(n+2\right)}-W_{m\left(n-3\right)}\right),$$
	which is just Theorem 2.3 in \cite{A14}.
\end{remark}

\begin{corollary}\label{corollary 3.4}
	For any positive integers $m,n,p$ and integers $q,l$, where $l\geq 1-m$  and $p^2+2q-2<p\sqrt{p^2+4q}$, we have
	$$\left(\sum_{k=n}^{\infty}\frac{1}{W_{mk+l}^4}\right)^{-1}\sim W_{mn+l}^4-W_{m\left(n-1\right)+l}^4+E+F+G,$$
	where
	$$E=4c_1^3c_2\left(-\alpha^2 q\right)^{mn+l}\left(1-\left(-\alpha^2 q\right)^{-m}-\frac{\left(\alpha^{4m}-1\right)^2}{\alpha^{3m}\left(\alpha^{5m}-\beta^m\right)}\right),$$
	$$F=c_1^2c_2^2q^{2\left(mn+l\right)}\left(-6+6q^{-2m}-\frac{10\left(\alpha^{4m}-1\right)^2}{\alpha^{2m}\left(\alpha^{6m}-\beta^{2m}\right)}+\frac{16\left(\alpha^{4m}-1\right)^3}{\alpha^{2m}\left(\alpha^{5m}-\beta^m\right)^2}\right),$$
	and
	{\footnotesize$$
		G=4c_1c_2^3\left(-\beta^2 q\right)^{mn+l}\left(1-\left(-\beta^2 q\right)^{-m}-\frac{5\left(\alpha^{4m}-1\right)^2}{\alpha^m \left(\alpha^{7m}-\beta^{3m}\right)}-\frac{16\left(\alpha^{4m}-1\right)^4}{\alpha^m \left(\alpha^{5m}-\beta^m\right)^3}+\frac{20\left(\alpha^{4m}-1\right)^3}{c_1^2\alpha^m\left(\alpha^{5m}-\beta^m\right)\left(\alpha^{6m}-\beta^{2m}\right)}\right).
		$$
	}
	%$$\begin{aligned}
		%G=&4c_1c_2^3\left(-\beta^2 q\right)^{mn+l}\left(1-\left(-\beta^2 q\right)^{-m}-\frac{5\left(\alpha^{4m}-1\right)^2}{\alpha^m \left(\alpha^{7m}-\beta^{3m}\right)}-\frac{16\left(\alpha^{4m}-1\right)^4}{\alpha^m \left(\alpha^{5m}-\beta^m\right)^3}\right.\\
		%&\left.+\frac{20\left(\alpha^{4m}-1\right)^3}{c_1^2\alpha^m\left(\alpha^{5m}-\beta^m\right)\left(\alpha^{6m}-\beta^{2m}\right)}\right).
		%\end{aligned}$$
	\end{corollary}
	{\bf Proof}.\ By Theorem \ref{theorem 3.1}, we have
	$$\begin{aligned}
		&\left(\sum_{k=n}^{\infty}\frac{1}{W_{mk+l}^4}\right)^{-1}\\
		=&W_{mn+l}^4-W_{m\left(n-1\right)+l}^4-\sum_{i=1}^{4}\binom{4}{i}c_1^{4-i}\left(-c_2\right)^{i}\left(\left(\alpha^{mn+l}\right)^{4-i}\left(\beta^{mn+l}\right)^i-\left(\alpha^{m\left(n-1\right)+l}\right)^{4-i}\left(\beta^{m\left(n-1\right)+l}\right)^i\right)\\
		&+\frac{c_1^4\alpha^{4\left(mn+l\right)}\left(\alpha^{4m}-1\right)}{\alpha^{4m}}\sum_{j=1}^{3}\left(-1\right)^j\left(\sum_{i=1}^{3}\binom{3+i}{3}\frac{c_2^{i}\beta^{i\left(mn+l\right)}\alpha^{mi}\left(\alpha^{4m}-1\right)}{c_1^{i}\alpha^{i\left(mn+l\right)}\left(\alpha^{m\left(4+i\right)}-\beta^{mi}\right)}\right)^j+O\left(\beta^{4\left(mn+l\right)}\right)\\
		=&W_{mn+l}^4-W_{m\left(n-1\right)+l}^4+E+F+G+H+O\left(\beta^{4\left(mn+l\right)}\right),
	\end{aligned}$$
	where
	$$E=4c_1^3c_2\left(-\alpha^2 q\right)^{mn+l}\left(1-\left(-\alpha^2 q\right)^{-m}-\frac{\left(\alpha^{4m}-1\right)^2}{\alpha^{3m}\left(\alpha^{5m}-\beta^m\right)}\right),$$
	$$F=c_1^2c_2^2q^{2\left(mn+l\right)}\left(-6+6q^{-2m}-\frac{10\left(\alpha^{4m}-1\right)^2}{\alpha^{2m}\left(\alpha^{6m}-\beta^{2m}\right)}+\frac{16\left(\alpha^{4m}-1\right)^3}{\alpha^{2m}\left(\alpha^{5m}-\beta^m\right)^2}\right),$$
	{\footnotesize$$
		G=4c_1c_2^3\left(-\beta^2 q\right)^{mn+l}\left(1-\left(-\beta^2 q\right)^{-m}-\frac{5\left(\alpha^{4m}-1\right)^2}{\alpha^m \left(\alpha^{7m}-\beta^{3m}\right)}-\frac{16\left(\alpha^{4m}-1\right)^4}{\alpha^m \left(\alpha^{5m}-\beta^m\right)^3}+\frac{20\left(\alpha^{4m}-1\right)^3}{c_1^2\alpha^m\left(\alpha^{5m}-\beta^m\right)\left(\alpha^{6m}-\beta^{2m}\right)}\right),
		$$
	}
	and
	$$\begin{aligned}
		H=&-c_2^4\beta^{4\left(mn+l\right)}
		+c_2^4\beta^{4\left(m\left(n-1\right)+l\right)}+\frac{100c_2^4\beta^{4\left(mn+l\right)}\left(\alpha^{4m}-1\right)^3}{\left(\alpha^{6m}-\beta^{2m}\right)^2}+\frac{400c_2^6\beta^{6\left(mn+l\right)}\alpha^{2m}\left(\alpha^{4m}-1\right)^3}{c_1^2\alpha^{2\left(mn+l\right)}\left(\alpha^{7m}-\beta^{3m}\right)^2}\\
		&+\frac{160c_2^4\beta^{4\left(mn+l\right)}\left(\alpha^{4m}-1\right)^3}{\left(\alpha^{m}-\beta^{m}\right)\left(\alpha^{7m}-\beta^{3m}\right)}+\frac{400c_2^5\beta^{5\left(mn+l\right)}\alpha^{m}\left(\alpha^{4m}-1\right)^3}{c_1\alpha^{mn+l}\left(\alpha^{6m}-\beta^{2m}\right)\left(\alpha^{7m}-\beta^{3m}\right)}\\
		&-\frac{1000c_2^6\beta^{6\left(mn+l\right)}\alpha^{2m}\left(\alpha^{4m}-1\right)^4}{c_1^2\alpha^{2\left(mn+l\right)}\left(\alpha^{6m}-\beta^{2m}\right)^3}-\frac{8000c_2^9\beta^{9\left(mn+l\right)}\alpha^{5m}\left(\alpha^{4m}-1\right)^4}{c_1^5\alpha^{5\left(mn+l\right)}\left(\alpha^{7m}-\beta^{3m}\right)^3}\\
		&-\frac{480c_2^4\beta^{4\left(mn+l\right)}\left(\alpha^{4m}-1\right)^4}{\left(\alpha^{5m}-\beta^{m}\right)^2\left(\alpha^{6m}-\beta^{2m}\right)}-\frac{960c_2^5\beta^{5\left(mn+l\right)}\alpha^{m}\left(\alpha^{4m}-1\right)^4}{c_1\alpha^{mn+l}\left(\alpha^{5m}-\beta^{m}\right)^2\left(\alpha^{7m}-\beta^{3m}\right)}\\
		&-\frac{1200c_2^5\beta^{5\left(mn+l\right)}\alpha^{m}\left(\alpha^{4m}-1\right)^4}{c_1\alpha^{mn+l}\left(\alpha^{5m}-\beta^{m}\right)\left(\alpha^{6m}-\beta^{2m}\right)^2}-\frac{6000c_2^7\beta^{7\left(mn+l\right)}\alpha^{3m}\left(\alpha^{4m}-1\right)^4}{c_1^3\alpha^{3\left(mn+l\right)}\left(\alpha^{6m}-\beta^{2m}\right)^2\left(\alpha^{7m}-\beta^{3m}\right)}\\
		&-\frac{4800c_2^7\beta^{7\left(mn+l\right)}\alpha^{3m}\left(\alpha^{4m}-1\right)^4}{c_1^3\alpha^{3\left(mn+l\right)}\left(\alpha^{5m}-\beta^{m}\right)\left(\alpha^{7m}-\beta^{3m}\right)^2}-\frac{12000c_2^8\beta^{8\left(mn+l\right)}\alpha^{4m}\left(\alpha^{4m}-1\right)^4}{c_1^4\alpha^{4\left(mn+l\right)}\left(\alpha^{6m}-\beta^{2m}\right)\left(\alpha^{7m}-\beta^{3m}\right)^2}\\
		&-\frac{4800c_2^6\beta^{6\left(mn+l\right)}\alpha^{2m}\left(\alpha^{4m}-1\right)^4}{c_1^2\alpha^{2\left(mn+l\right)}\left(\alpha^{5m}-\beta^{m}\right)\left(\alpha^{6m}-\beta^{2m}\right)\left(\alpha^{7m}-\beta^{3m}\right)}.
	\end{aligned}$$
	Note that $p>0$, it's easy to prove that $|\alpha|>|\beta|$, furthermore, by  $p^2+2q-2<p\sqrt{p^2+4q}$, it's easy to prove that $|\beta|<1$. Therefore, for any positive integers $m$ and $n$ with $n<m$, we have $\left|\alpha^n\right|>\left|\beta^n\right|>\left|\beta^m\right|$, thus
	$$
	\begin{aligned}
		&\lim_{n\rightarrow \infty}
		\left(\left(\sum_{k=n}^{\infty}\frac{1}{W_{mk+l}^4}\right)^{-1}-\left(W_{mn+l}^4-W_{m\left(n-1\right)+l}^4+E+F+G\right)\right)\\
		=&\lim_{n\rightarrow \infty}
		\left(H+O\left(\beta^{4\left(mn+l\right)}\right)\right)\\
		=&0.
	\end{aligned}
	$$
	Therefore, by Definition \ref{definition 2.1}, we can get
	$$\left(\sum_{k=n}^{\infty}\frac{1}{W_{mk+l}^4}\right)^{-1}\sim W_{mn+l}^4-W_{m\left(n-1\right)+l}^4+E+F+G,$$
	where
	$$E=4c_1^3c_2\left(-\alpha^2 q\right)^{mn+l}\left(1-\left(-\alpha^2 q\right)^{-m}-\frac{\left(\alpha^{4m}-1\right)^2}{\alpha^{3m}\left(\alpha^{5m}-\beta^m\right)}\right),$$
	$$F=c_1^2c_2^2q^{2\left(mn+l\right)}\left(-6+6q^{-2m}-\frac{10\left(\alpha^{4m}-1\right)^2}{\alpha^{2m}\left(\alpha^{6m}-\beta^{2m}\right)}+\frac{16\left(\alpha^{4m}-1\right)^3}{\alpha^{2m}\left(\alpha^{5m}-\beta^m\right)^2}\right)$$
	and
	{\footnotesize$$
		G=4c_1c_2^3\left(-\beta^2 q\right)^{mn+l}\left(1-\left(-\beta^2 q\right)^{-m}-\frac{5\left(\alpha^{4m}-1\right)^2}{\alpha^m \left(\alpha^{7m}-\beta^{3m}\right)}-\frac{16\left(\alpha^{4m}-1\right)^4}{\alpha^m \left(\alpha^{5m}-\beta^m\right)^3}+\frac{20\left(\alpha^{4m}-1\right)^3}{c_1^2\alpha^m\left(\alpha^{5m}-\beta^m\right)\left(\alpha^{6m}-\beta^{2m}\right)}\right).
		$$
	}

	This completes the proof of Corollary \ref{corollary 3.4}.
	
	$\hfill\Box$
	
	The following Remark 3.4 shows that the corresponding result is just Theorem 2.4 in \cite{A14} when $a=0,b=1,p=A\in \mathbb{Z}^+,q=\pm 1$ and $l=0$ in Corollary \ref{corollary 3.4}.
	%In the following Remark 3.3, it's explained that Corollary \ref{corollary 3.4} is just Theorem 2.4 in \cite{A14} when $a=0,b=1,p=A\in \mathbb{Z}^+,q=\pm 1$ and $l=0$.
	\begin{remark}
		By taking $a=0,b=1,p=A\in \mathbb{Z}^+,q=B=\pm 1$ and $l=0$ in Corollary \ref{corollary 3.4}, we have $c_1=c_2=\frac{1}{\alpha-\beta}$, and then 
		$$\begin{aligned}
			E&=\frac{4B^{mn}\alpha^{2mn}}{\left(\alpha-\beta\right)^4}\left(1-B^{-m}\alpha^{-2m}-\frac{\left(\alpha^{4m}-1\right)^2}{\alpha^{3m}\left(\alpha^{5m}-\beta^m\right)}\right)\\
			&=\frac{4B^{mn}\alpha^{2mn}}{\left(\alpha-\beta\right)^4}\left(\frac{\left(1-B^{m}\alpha^{-2m}\right)\alpha^{2m}\left(\alpha^{6m}-\alpha^m\beta^m\right)-\left(\alpha^{4m}-1\right)^2}{\alpha^{2m}\left(\alpha^{6m}-\alpha^m\beta^m\right)}\right),
		\end{aligned}$$
		$$\begin{aligned}
			F&=\frac{1}{\left(\alpha-\beta\right)^4}\left(-\frac{10\left(\alpha^{4m}-1\right)^2}{\alpha^{2m}\left(\alpha^{6m}-\beta^{2m}\right)}+\frac{16\left(\alpha^{4m}-1\right)^3}{\alpha^{2m}\left(\alpha^{5m}-\beta^{m}\right)^2}\right)\\
			&=\frac{\left(\alpha^{4m}-1\right)^2}{\left(\alpha-\beta\right)^4}\left(\frac{16\left(\alpha^{4m}-1\right)}{\left(\alpha^{6m}-\alpha^m\beta^{m}\right)^2}-\frac{10}{\alpha^{8m}-\alpha^{2m}\beta^{2m}}\right)
		\end{aligned}$$
		and
		{\small$$
			G=\frac{4\beta^{mn}B^{mn}}{\left(\alpha-\beta^3\right)}\left(1-\beta^{-2m}B^{-m}-\frac{5\left(\alpha^{4m}-1\right)^2}{\alpha^m \left(\alpha^{7m}-\beta^{3m}\right)}-\frac{16\left(\alpha^{4m}-1\right)^4}{\alpha^m \left(\alpha^{5m}-\beta^m\right)^3}+\frac{20\left(\alpha^{4m}-1\right)^3}{c_1^2\alpha^m\left(\alpha^{5m}-\beta^m\right)\left(\alpha^{6m}-\beta^{2m}\right)}\right).
			$$}
		%$$\begin{aligned}
			%G=&\frac{4\beta^{mn}B^{mn}}{\left(\alpha-\beta^3\right)}\left(1-\beta^{-2m}B^{-m}-\frac{5\left(\alpha^{4m}-1\right)^2}{\alpha^m \left(\alpha^{7m}-\beta^{3m}\right)}-\frac{16\left(\alpha^{4m}-1\right)^4}{\alpha^m \left(\alpha^{5m}-\beta^m\right)^3}\right.\\
			%&\left.+\frac{20\left(\alpha^{4m}-1\right)^3}{c_1^2\alpha^m\left(\alpha^{5m}-\beta^m\right)\left(\alpha^{6m}-\beta^{2m}\right)}\right).
			%\end{aligned}$$
			Note that $\alpha\beta=B$, it's easy to know that $\alpha^{m}
			\beta^m=B^m$ and $\alpha^{2m}
			\beta^{2m}=B^{2m}=1$, thus,
			$$
			E=\frac{4B^{mn}\alpha^{2mn}}{\left(\alpha-\beta\right)^4}\left(\frac{\left(1-B^{m}\beta^{2m}\right)\alpha^{2m}\left(\alpha^{6m}-B^m\right)-\left(\alpha^{4m}-1\right)^2}{\alpha^{2m}\left(\alpha^{6m}-B^m\right)}\right)
			$$
			and
			$$\begin{aligned}
				F=\frac{\left(\alpha^{4m}-1\right)^2}{\left(\alpha-\beta\right)^4}\left(\frac{16\left(\alpha^{4m}-1\right)}{\left(\alpha^{6m}-B^m\right)^2}-\frac{10}{\alpha^{8m}-1}\right).
			\end{aligned}$$
			Therefore, E and F are just $\delta$ and $V_m$ in the proof process of Theorem 2.4 in \cite{A14},  respectively, furthermore,
			by the proof of Theorem 2.4 in \cite{A14}, we have  
			$$E=4B^{mn}U_m\left(W_{m\left(n+1\right)}^2-B^mW_{\left(m\left(n-2\right)\right)}^2\right)+O\left(\frac{1}{\alpha^{2m\left(n-2\right)}}\right),$$
			where $U_m=\frac{W_m^2}{\left(1-B^m\alpha^{6m}\right)\left(1-B^m\beta^{6m}\right)}$.
			By $\left|\beta\right|<1$ and $\left|\alpha\right|=\left|B\beta^{-1}\right|=\left|\beta^{-1}\right|>1$, it's easy to get
			$$
			\begin{aligned}
				&\lim_{n\rightarrow \infty}
				\left(\left(\sum_{k=n}^{\infty}\frac{1}{W_{mk}^4}\right)^{-1}-\left(W_{mn}^4-W_{m\left(n-1\right)}^4+4B^{mn}U_m\left(W_{m\left(n+1\right)}^2-B^mW_{m\left(n-2\right)}^2\right)+V_m\right)\right)\\
				=&\lim_{n\rightarrow \infty}
				\left(O\left(\frac{1}{\alpha^{2m\left(n-2\right)}}\right)+G\right)\\
				=&0.
			\end{aligned}
			$$
			Therefore, by Definition \ref{definition 2.1}, we can get
			$$\left(\sum_{k=n}^{\infty}\frac{1}{W_{mk}^4}\right)^{-1}\sim W_{mn}^4-W_{m\left(n-1\right)}^4+4B^{mn}U_m\left(W_{m\left(n+1\right)}^2-B^mW_{m\left(n-2\right)}^2\right)+V_m,$$
			which is just Theorem 2.4 in \cite{A14}.
		\end{remark}

		\section{The asymptotic estimation value of $\left(\sum\limits_{k=n}^{\infty}\frac{\left(-1\right)^k}{W_{mk+l}^d}\right)^{-1}$}
		In this section, for the generalized Fibonacci sequence $\left\{W_n\left(a,b,p,q\right)\right\}$, we establish the corresponding connection of $\left(\sum\limits_{k=n}^{\infty}\frac{\left(-1\right)^k}{W_{mk+l}^d}\right)^{-1}$ and its asymptotic estimation(Theorem \ref{theorem 4.1}), and then obtain the more precise asymptotic estimation values of $\left(\sum\limits_{k=n}^{\infty}\frac{\left(-1\right)^k}{W_{mk+l}^d}\right)^{-1}$ for $d=1,2,3,4$, respectively(Corollaries 4.1–4.4).
		\subsection{The connection of $\left(\sum\limits_{k=n}^{\infty}\frac{\left(-1\right)^k}{W_{mk+l}^d}\right)^{-1}$ and its asymptotic estimation}
		In this subsection, the connection of $\left(\sum\limits_{k=n}^{\infty}\frac{\left(-1\right)^k}{W_{mk+l}^d}\right)^{-1}$ and its asymptotic estimation is given as the following 
		\begin{theorem}\label{theorem 4.1}
			For any positive integers $m,n,d,p$ and integers $q,l$, where $l\geq 1-m$  and $p^2+2q-2<p\sqrt{p^2+4q}$, we have
			$$
			\small \begin{aligned}
				&\left(\sum_{k=n}^{\infty}\frac{\left(-1\right)^k}{W_{mk+l}^d}\right)^{-1}\\
				=&\left(-1\right)^n\left(W_{mn+l}^d+W_{m\left(n-1\right)+l}^d\right)\\
				&-\left(-1\right)^n\sum_{i=1}^{d}\binom{d}{i}c_1^{d-i}\left(-c_2\right)^{i}\left(\left(\alpha^{mn+l}\right)^{d-i}\left(\beta^{mn+l}\right)^i+\left(\alpha^{m\left(n-1\right)+l}\right)^{d-i}\left(\beta^{m\left(n-1\right)+l}\right)^i\right)\\
				&+\frac{\left(-1\right)^{n}c_1^d\alpha^{d\left(mn+l\right)}\left(\alpha^{md}+1\right)}{\alpha^{md}}\sum_{j=1}^{d-1}\left(-1\right)^j\left(\sum_{i=1}^{d-1}\binom{d-1+i}{d-1}\frac{c_2^{i}\beta^{i\left(mn+l\right)}\alpha^{mi}\left(\alpha^{md}+1\right)}{c_1^{i}\alpha^{i\left(mn+l\right)}\left(\alpha^{m\left(d+i\right)}+\beta^{mi}\right)}\right)^j+O\left(\beta^{d\left(mn+l\right)}\right).
			\end{aligned}
			$$
		\end{theorem}
		{\bf Proof}.\ By Lemma \ref{lemma 2.3}, we have
		$$
		\begin{aligned}
			&\sum_{k=n}^{\infty}\frac{\left(-1\right)^k}{W_{mk+l}^d\left(a,b,p,q\right)}\\
			=&\frac{1}{c_1^d}\left(\sum_{k=n}^{\infty}\frac{\left(-1\right)^k}{\alpha^{d\left(mk+l\right)}}+\sum_{i=1}^{d-1}\binom{d-1+i}{d-1}\left(\frac{c_2}{c_1}\right)^i\sum_{k=n}^{\infty}\frac{\left(-1\right)^k\beta^{i\left(mk+l\right)}}{\alpha^{\left(mk+l\right)\left(d+i\right)}}+\sum_{k=n}^{\infty}\left(-1\right)^{k}O\left(\frac{\beta^{d\left(mk+l\right)}}{\alpha^{2d\left(mk+l\right)}}\right)\right).
		\end{aligned}
		$$
		Note that
		$$\begin{aligned}
			\sum_{k=n}^{\infty}\frac{\left(-1\right)^k}{\alpha^{d\left(mk+l\right)}}&
			=\frac{1}{\alpha^{dl}}\sum_{k=n}^{\infty}\left(\frac{-1}{\alpha^{dm}}\right)^k\\
			&=\frac{1}{\alpha^{dl}}\left(\frac{\left(-1\right)^n}{\alpha^{dmn}}+\frac{\left(-1\right)^{n+1}}{\alpha^{dmn+dm}}+\frac{\left(-1\right)^{n+2}}{\alpha^{dmn+2dm}}+\cdots\right)\\
			&=\frac{\left(-1\right)^n}{\alpha^{dmn+dl}}\left(1+\frac{-1}{\alpha^{dm}}+\left(\frac{-1}{\alpha^{dm}}\right)^2+\cdots\right)\\
			&=\frac{\left(-1\right)^n}{\alpha^{dmn+dl}}\cdot\frac{1}{1-\frac{-1}{\alpha^{dm}}}\\
			&=\frac{\left(-1\right)^n\alpha^{md}}{\alpha^{d\left(mn+l\right)}\left(\alpha^{md}+1\right)},
		\end{aligned}$$
		$$\begin{aligned}
			\sum_{k=n}^{\infty}\frac{\left(-1\right)^k\beta^{i\left(mk+l\right)}}{\alpha^{\left(mk+l\right)\left(d+i\right)}}&=\frac{\beta^{li}}{\alpha^{l\left(d+i\right)}}\sum_{k=n}^{\infty}\left(\frac{-\beta^{mi}}{\alpha^{m\left(d+i\right)}}\right)^k\\
			&=\frac{\beta^{li}}{\alpha^{l\left(d+i\right)}}\left(\frac{\left(-1\right)^n\beta^{mni}}{\alpha^{mn\left(d+i\right)}}+\frac{\left(-1\right)^{n+1}\beta^{mni+mi}}{\alpha^{mn\left(d+i\right)+m\left(d+i\right)}}+\frac{\left(-1\right)^{n+2}\beta^{mni+2mi}}{\alpha^{mn\left(d+i\right)+2m\left(d+i\right)}}+\cdots\right)\\
			&=\frac{\beta^{li}}{\alpha^{l\left(d+i\right)}}\cdot\frac{\left(-1\right)^n\beta^{mni}}{\alpha^{mn\left(d+i\right)}}\left(1+\frac{-\beta^{mi}}{\alpha^{m\left(d+i\right)}}+\left(\frac{-\beta^{mi}}{\alpha^{m\left(d+i\right)}}\right)^2+\cdots\right)\\
			&=\frac{\left(-1\right)^n\beta^{i\left(mn+l\right)}}{\alpha^{\left(mn+l\right)\left(d+i\right)}}\cdot\frac{1}{1-\frac{-\beta^{mi}}{\alpha^{m\left(d+i\right)}}}\\
			&=\frac{\left(-1\right)^n\beta^{i\left(mn+l\right)}\alpha^{m\left(d+i\right)}}{\alpha^{\left(mn+l\right)\left(d+i\right)}\left(\alpha^{m\left(d+i\right)}+\beta^{mi}\right)}
		\end{aligned}$$
		and by Lemma \ref{lemma 2.2}, it's easy to get
		$$\sum_{k=n}^{\infty}\left(-1\right)^kO\left(\frac{\beta^{d\left(mk+l\right)}}{\alpha^{2d\left(mk+l\right)}}\right)=\sum_{k=n}^{\infty}O\left(\frac{\beta^{d\left(mk+l\right)}}{\alpha^{2d\left(mk+l\right)}}\right)=O\left(\frac{\beta^{d\left(mn+l\right)}}{\alpha^{2d\left(mn+l\right)}}\right),$$
		thus
		$$\begin{aligned}
			&\sum_{k=n}^{\infty}\frac{\left(-1\right)^k}{W_{mk+l}^d}\\
			=&\frac{1}{c_1^d}\left(\frac{\left(-1\right)^n\alpha^{md}}{\alpha^{d\left(mn+l\right)}\left(\alpha^{md}+1\right)}+\sum_{i=1}^{d-1}\binom{d-1+i}{d-1}\frac{\left(-1\right)^nc_2^i\beta^{i\left(mn+l\right)}\alpha^{m\left(d+i\right)}}{c_1^i\alpha^{\left(mn+l\right)\left(d+i\right)}\left(\alpha^{m\left(d+i\right)}+\beta^{mi}\right)}+O\left(\frac{\beta^{d\left(mn+l\right)}}{\alpha^{2d\left(mn+l\right)}}\right)\right)\\
			=&\frac{\left(-1\right)^n\alpha^{md}}{c_1^d\alpha^{d\left(mn+l\right)}\left(\alpha^{md}+1\right)}\left(1+\sum_{i=1}^{d-1}\binom{d-1+i}{d-1}\frac{c_2^i\beta^{i\left(mn+l\right)}\alpha^{mi}\left(\alpha^{md}+1\right)}{c_1^i\alpha^{i\left(mn+l\right)}\left(\alpha^{m\left(d+i\right)}+\beta^{mi}\right)}+O\left(\frac{\beta^{d\left(mn+l\right)}}{\alpha^{d\left(mn+l\right)}}\right)\right),
		\end{aligned}$$
		and then
		$$\begin{aligned}
			&\left(\sum_{k=n}^{\infty}\frac{\left(-1\right)^k}{W_{mk+l}^d}\right)^{-1}\\
			=&\frac{\left(-1\right)^nc_1^d\alpha^{d\left(mn+l\right)}\left(\alpha^{md}+1\right)}{\alpha^{md}}\left(1+\sum_{i=1}^{d-1}\binom{d-1+i}{d-1}\frac{c_2^i\beta^{i\left(mn+l\right)}\alpha^{mi}\left(\alpha^{md}+1\right)}{c_1^i\alpha^{i\left(mn+l\right)}\left(\alpha^{m\left(d+i\right)}+\beta^{mi}\right)}+O\left(\frac{\beta^{d\left(mn+l\right)}}{\alpha^{d\left(mn+l\right)}}\right)\right)^{-1}.
		\end{aligned}$$
		And by $p>0$, it's easy to prove $\left|\frac{\beta}{\alpha}\right|<1$. Thus, for any positive integer $n$, we have
		$$\lim_{n\rightarrow \infty}\left(\sum_{i=1}^{d-1}\binom{d-1+i}{d-1}\frac{c_2^i\beta^{il}\alpha^{mi}\left(\alpha^{md}+1\right)}{c_1^i\alpha^{il}\left(\alpha^{m\left(d+i\right)}+\beta^{mi}\right)}\cdot\left(\frac{\beta}{\alpha}\right)^{imn}+O\left(\frac{\beta^{d\left(mn+l\right)}}{\alpha^{d\left(mn+l\right)}}\right)\right)=0,$$
		then by the definition, it's easy to know that there exists a sufficiently large positive integer $n$ such that
		$$\left|\sum_{i=1}^{d-1}\binom{d-1+i}{d-1}\frac{c_2^i\beta^{i\left(mn+l\right)}\alpha^{mi}\left(\alpha^{md}+1\right)}{c_1^{i}\alpha^{i\left(mn+l\right)}\left(\alpha^{m\left(d+i\right)}+\beta^{mi}\right)}+O\left(\frac{\beta^{d\left(mn+l\right)}}{\alpha^{d\left(mn+l\right)}}\right)\right|<1.$$
		Furthermore, by Lemma \ref{definition 2.1}, we have
		$$\small \begin{aligned}
			&\left(\sum_{k=n}^{\infty}\frac{\left(-1\right)^k}{W_{mk+l}^d}\right)^{-1}\\
			=&\frac{\left(-1\right)^nc_1^d\alpha^{d\left(mn+l\right)}\left(\alpha^{md}+1\right)}{\alpha^{md}}\left(1+\sum_{j=1}^{d-1}\left(-1\right)^j\left(\sum_{i=1}^{d-1}\binom{d-1+i}{d-1}\frac{c_2^{i}\beta^{i\left(mn+l\right)}\alpha^{mi}\left(\alpha^{md}+1\right)}{c_1^{i}\alpha^{i\left(mn+l\right)}\left(\alpha^{m\left(d+i\right)}+\beta^{mi}\right)}\right)^j+O\left(\frac{\beta^{d\left(mn+l\right)}}{\alpha^{d\left(mn+l\right)}}\right) \right)\\
			=&\left(-1\right)^n\left(c_1^d\alpha^{dmn+dl}+c_1\alpha^{dm\left(n-1\right)+dl}\right)\\
			&+\frac{\left(-1\right)^nc_1^d\alpha^{d\left(mn+l\right)}\left(\alpha^{md}+1\right)}{\alpha^{md}}\sum_{j=1}^{d-1}\left(-1\right)^j\left(\sum_{i=1}^{d-1}\binom{d-1+i}{d-1}\frac{c_2^{i}\beta^{i\left(mn+l\right)}\alpha^{mi}\left(\alpha^{md}+1\right)}{c_1^{i}\alpha^{i\left(mn+l\right)}\left(\alpha^{m\left(d+i\right)}+\beta^{mi}\right)}\right)^j+O\left(\beta^{d\left(mn+l\right)}\right).
		\end{aligned}$$
		Note that for any positive integer $e$, we have
		$$W_{me+l}^d=c_1^d\alpha^{dme+dl}+\sum_{i=1}^{d}\binom{d}{i}c_1^{d-i}\left(-c_2\right)^{i}\left(\alpha^{me+l}\right)^{d-i}\left(\beta^{me+l}\right)^i,$$
		thus
		$$\small \begin{aligned}
			&\left(\sum_{k=n}^{\infty}\frac{\left(-1\right)^k}{W_{mk+l}^d}\right)^{-1}\\
			=&\left(-1\right)^n\left(W_{mn+l}^d+W_{m\left(n-1\right)+l}^d\right)\\
			&-\left(-1\right)^n\sum_{i=1}^{d}\binom{d}{i}c_1^{d-i}\left(-c_2\right)^{i}\left(\left(\alpha^{mn+l}\right)^{d-i}\left(\beta^{mn+l}\right)^i+\left(\alpha^{m\left(n-1\right)+l}\right)^{d-i}\left(\beta^{m\left(n-1\right)+l}\right)^i\right)\\
			&+\frac{\left(-1\right)^{n}c_1^d\alpha^{d\left(mn+l\right)}\left(\alpha^{md}+1\right)}{\alpha^{md}}\sum_{j=1}^{d-1}\left(-1\right)^j\left(\sum_{i=1}^{d-1}\binom{d-1+i}{d-1}\frac{c_2^{i}\beta^{i\left(mn+l\right)}\alpha^{mi}\left(\alpha^{md}+1\right)}{c_1^{i}\alpha^{i\left(mn+l\right)}\left(\alpha^{m\left(d+i\right)}+\beta^{mi}\right)}\right)^j+O\left(\beta^{d\left(mn+l\right)}\right).
		\end{aligned}$$
		
		This completes the proof of Theorem \ref{theorem 4.1}.
		
		$\hfill\Box$

		\subsection{Some Corollaries}
		In this subsection, for some special $d$, we can obtain more precise asymptotic estimate values of $\left(\sum\limits_{k=n}^{\infty}\frac{\left(-1\right)^k}{W_{mk+l}^d}\right)^{-1}$ as the following Corollaries \ref{corollary 4.1}-\ref{corollary 4.4}.
		
		\begin{corollary}\label{corollary 4.1}
			For any positive integers $m,n,p$ and integers $q,l$, where $l\geq 1-m$  and $p^2+2q-2<p\sqrt{p^2+4q}$, we have
			$$\left(\sum_{k=n}^{\infty}\frac{\left(-1\right)^k}{W_{mk+l}}\right)^{-1}\sim \left(-1\right)^n\left(W_{mn+l}+W_{m\left(n-1\right)+l}\right).$$
		\end{corollary}
		{\bf Proof}.\ By Theorem \ref{theorem 4.1}, we have
		$$\small 
		\left(\sum_{k=n}^{\infty}\frac{\left(-1\right)^k}{W_{mk+l}}\right)^{-1}=\left(-1\right)^n\left(W_{mn+l}+W_{m\left(n-1\right)+l}\right)-\left(-1\right)^nc_2\beta^{m\left(n-1\right)+l}\left(\beta^m+1\right)+O\left(\beta^{mn+l}\right),
		$$
		and by $p^2+2q-2<p\sqrt{p^2+4q}$, it's easy to prove that $|\beta|<1$, hence,
		$$
		\begin{aligned}
			&\lim_{n\rightarrow \infty}
			\left(\left(\sum_{k=n}^{\infty}\frac{\left(-1\right)^k}{W_{mk+l}}\right)^{-1}-\left(-1\right)^n\left(W_{mn+l}+W_{m\left(n-1\right)+l}\right)\right)\\
			=&\lim_{n\rightarrow \infty}
			\left(-\left(-1\right)^nc_2\beta^{m\left(n-1\right)+l}\left(\beta^m+1\right)+O\left(\beta^{mn+l}\right)\right)\\
			=&0.
		\end{aligned}
		$$
		Therefore, by Definition \ref{definition 2.1}, we can get
		$$\left(\sum_{k=n}^{\infty}\frac{\left(-1\right)^k}{W_{mk+l}}\right)^{-1}\sim \left(-1\right)^n\left(W_{mn+l}+W_{m\left(n-1\right)+l}\right).$$
		
		This completes the proof of Corollary \ref{corollary 4.1}.
		
		$\hfill\Box$

		\begin{corollary}\label{corollary 4.2}
			For any positive integers $m,n,p$ and integers $q,l$, where $l\geq 1-m$  and $p^2+2q-2<p\sqrt{p^2+4q}$, we have
			$$\left(\sum_{k=n}^{\infty}\frac{\left(-1\right)^k}{W_{mk+l}^2}\right)^{-1}\sim \left(-1\right)^n\left(W_{mn+l}^2+W_{m\left(n-1\right)+l}^2+H\right),$$
			where $H=2c_1c_2\left(-q\right)^{mn+l}\left(1+\left(-q\right)^{m}-\frac{\left(\alpha^{2m}+1\right)^2}{\alpha^m\left(\alpha^{3m}+\beta^m\right)}\right)$.
			
		\end{corollary}
		{\bf Proof}.\ By Theorem \ref{theorem 4.1}, we have
		$$\begin{aligned}
			&\left(\sum_{k=n}^{\infty}\frac{\left(-1\right)^k}{W_{mk+l}^2}\right)^{-1}\\
			=&\left(-1\right)^n\left(W_{mn+l}^2+W_{m\left(n-1\right)+l}^2\right)\\
			&-\left(-1\right)^n\sum_{i=1}^{2}\binom{2}{i}c_1^{2-i}\left(-c_2\right)^{i}\left(\left(\alpha^{mn+l}\right)^{2-i}\left(\beta^{mn+l}\right)^i+\left(\alpha^{m\left(n-1\right)+l}\right)^{2-i}\left(\beta^{m\left(n-1\right)+l}\right)^i\right)\\
			&-\frac{\left(-1\right)^{n}c_1^2\alpha^{2\left(mn+l\right)}\left(\alpha^{2m}+1\right)}{\alpha^{2m}}\cdot \frac{2c_2\beta^{mn+l}\alpha^{m}\left(\alpha^{2m}+1\right)}{c_1\alpha^{mn+l}\left(\alpha^{3m}+\beta^{m}\right)}+O\left(\beta^{2\left(mn+l\right)}\right)\\
			=&\left(-1\right)^n\left(W_{mn+l}^2+W_{m\left(n-1\right)+l}^2+H\right)-\left(-1\right)^nc_2^2\beta^{2\left(m\left(n-1\right)+l\right)}\left(\beta^{2m}+1\right)+O\left(\beta^{2\left(mn+l\right)}\right),
		\end{aligned}$$
		where $H=2c_1c_2\left(-q\right)^{mn+l}\left(1+\left(-q\right)^{m}-\frac{\left(\alpha^{2m}+1\right)^2}{\alpha^m\left(\alpha^{3m}+\beta^m\right)}\right)$.
		By $p^2+2q-2<p\sqrt{p^2+4q}$, it's easy to prove that $|\beta|<1$, hence,
		$$
		\begin{aligned}
			&\lim_{n\rightarrow \infty}
			\left(\left(\sum_{k=n}^{\infty}\frac{\left(-1\right)^k}{W_{mk+l}^2}\right)^{-1}-\left(-1\right)^n\left(W_{mn+l}^2+W_{m\left(n-1\right)+l}^2+H\right)\right)\\
			=&\lim_{n\rightarrow \infty}
			\left(-\left(-1\right)^nc_2^2\beta^{2\left(m\left(n-1\right)+l\right)}\left(\beta^{2m}+1\right)+O\left(\beta^{2\left(mn+l\right)}\right)\right)\\
			=&0.
		\end{aligned}
		$$
		Therefore, by Definition \ref{definition 2.1}, we have
		$$\left(\sum_{k=n}^{\infty}\frac{\left(-1\right)^k}{W_{mk+l}^2}\right)^{-1}\sim \left(-1\right)^n\left(W_{mn+l}^2+W_{m\left(n-1\right)+l}^2+H\right),$$
		where $H=2c_1c_2\left(-q\right)^{mn+l}\left(1+\left(-q\right)^{m}-\frac{\left(\alpha^{2m}+1\right)^2}{\alpha^m\left(\alpha^{3m}+\beta^m\right)}\right)$.
		
		This completes the proof of Corollary \ref{corollary 4.2}.
		
		$\hfill\Box$
		
		\begin{corollary}\label{corollary 4.3}
			For any positive integers $m,n,p$ and integers $q,l$, where $l\geq 1-m$  and $p^2+2q-2<p\sqrt{p^2+4q}$, we have
			$$\left(\sum_{k=n}^{\infty}\frac{\left(-1\right)^k}{W_{mk+l}^3}\right)^{-1}\sim \left(-1\right)^n\left(W_{mn+l}^3+W_{m\left(n-1\right)+l}^3+I+J\right),$$
			where
			$$I=3c_1^2c_2\left(-\alpha q\right)^{mn+l}\left(1+\left(-\alpha q\right)^{-m}-\frac{\left(\alpha^{3m}+1\right)^2}{\alpha^{2m}\left(\alpha^{4m}+\beta^m\right)}\right)
			$$
			and 
			$$J=3c_1c_2^2\left(-\beta q\right)^{mn+l}\left(-1-\left(-\beta q\right)^{-m}-\frac{2\left(\alpha^{3m}+1\right)^2}{\alpha^m\left(\alpha^{5m}+\beta^m\right)}+\frac{3\left(\alpha^{3m}+1\right)^3}{\alpha^m\left(\alpha^{4m}+\beta^m\right)^2}\right).$$
		\end{corollary}
		{\bf Proof}.\ By Theorem \ref{theorem 4.1}, we have
		$$\small \begin{aligned}
			&\left(\sum_{k=n}^{\infty}\frac{\left(-1\right)^k}{W_{mk+l}^3}\right)^{-1}\\
			=&\left(-1\right)^n\left(W_{mn+l}^3+W_{m\left(n-1\right)+l}^3\right)\\
			&-\left(-1\right)^n\sum_{i=1}^{3}\binom{3}{i}c_1^{3-i}\left(-c_2\right)^{i}\left(\left(\alpha^{mn+l}\right)^{3-i}\left(\beta^{mn+l}\right)^i+\left(\alpha^{m\left(n-1\right)+l}\right)^{3-i}\left(\beta^{m\left(n-1\right)+l}\right)^i\right)\\
			&+\frac{\left(-1\right)^{n}c_1^3\alpha^{3\left(mn+l\right)}\left(\alpha^{3m}+1\right)}{\alpha^{3m}}\sum_{j=1}^{2}\left(-1\right)^j\left(\sum_{i=1}^{2}\binom{2+i}{2}\frac{c_2^{i}\beta^{i\left(mn+l\right)}\alpha^{mi}\left(\alpha^{3m}+1\right)}{c_1^{i}\alpha^{i\left(mn+l\right)}\left(\alpha^{m\left(3+i\right)}+\beta^{mi}\right)}\right)^j+O\left(\beta^{3\left(mn+l\right)}\right)\\
			=&\left(-1\right)^n\left(W_{mn+l}^3+W_{m\left(n-1\right)+l}^3+I+J\right)+\left(-1\right)^nc_2^3\beta^{3\left(m\left(n-1\right)+l\right)}\left(\beta^{3m}+1+\frac{36\beta^{3m}\alpha^m\left(\alpha^{3m}+1\right)^3}{\left(\alpha^{4m}+\beta^{m}\right)\left(\alpha^{5m}+\beta^{2m}\right)}\right)\\
			&+\frac{36\left(-1\right)^nc_2^4\beta^{4\left(mn+l\right)}\alpha^m\left(\alpha^{3m}+1\right)^3}{c_1\alpha^{mn+l}\left(\alpha^{5m}+\beta^{2m}\right)^2}+O\left(\beta^{3\left(mn+l\right)}\right),
		\end{aligned}$$
		where
		$$I=3c_1^2c_2\left(-\alpha q\right)^{mn+l}\left(1+\left(-\alpha q\right)^{-m}-\frac{\left(\alpha^{3m}+1\right)^2}{\alpha^{2m}\left(\alpha^{4m}+\beta^m\right)}\right)$$
		and
		$$J=3c_1c_2^2\left(-\beta q\right)^{mn+l}\left(-1-\left(-\beta q\right)^{-m}-\frac{2\left(\alpha^{3m}+1\right)^2}{\alpha^m\left(\alpha^{5m}+\beta^m\right)}+\frac{3\left(\alpha^{3m}+1\right)^3}{\alpha^m\left(\alpha^{4m}+\beta^m\right)^2}\right).$$
		By $p>0$, it's easy to prove that $|\alpha|>|\beta|$. And by  $p^2+2q-2<p\sqrt{p^2+4q}$, it's easy to prove that $|\beta|<1$. Then, $|\alpha|>|\beta|>|\beta^4|$.  Furthermore,
		{\scriptsize $$
			\begin{aligned}
				&\lim_{n\rightarrow \infty}
				\left(\left(\sum_{k=n}^{\infty}\frac{\left(-1\right)^k}{W_{mk+l}^3}\right)^{-1}-\left(-1\right)^n\left(W_{mn+l}^3+W_{m\left(n-1\right)+l}^3+I+J\right)\right)\\
				%=&\lim_{n\rightarrow \infty} \left(\left(-1\right)^nc_2^3\beta^{3\left(mn+l\right)}+\left(-1\right)^nc_2^3\beta^{3\left(m\left(n-1\right)+l\right)}+\frac{36\left(-1\right)^nc_2^3\beta^{3\left(mn+l\right)}\alpha^m\left(\alpha^{3m}-1\right)^3}{\left(\alpha^{4m}-\beta^{m}\right)\left(\alpha^{5m}-\beta^{2m}\right)}\right.\\
				%&\left.+\frac{36\left(-1\right)^nc_2^4\beta^{4\left(mn+l\right)}\alpha^m\left(\alpha^{3m}-1\right)^3}{c_1\alpha^{mn+l}\left(\alpha^{5m}-\beta^{2m}\right)^2} +O\left(\beta^{3\left(mn+l\right)}\right)\right)\\
				=&\lim_{n\rightarrow \infty}
				\left(\left(-1\right)^nc_2^3\beta^{3\left(m\left(n-1\right)+l\right)}\left(\beta^{3m}+1+\frac{36\beta^{3m}\alpha^m\left(\alpha^{3m}+1\right)^3}{\left(\alpha^{4m}+\beta^{m}\right)\left(\alpha^{5m}+\beta^{2m}\right)}\right)+\left(\frac{\beta^4}{\alpha}\right)^{mn+l}\cdot \frac{36\left(-1\right)^nc_2^4\alpha^m\left(\alpha^{3m}+1\right)^3}{c_1\left(\alpha^{5m}+\beta^{2m}\right)^2}+O\left(\beta^{3\left(mn+l\right)}\right) \right)\\
				=&0.
			\end{aligned}
			$$}
		%$$
		%\begin{aligned}
		%&\lim_{n\rightarrow \infty}
		%\left(\left(\sum_{k=n}^{\infty}\frac{\left(-1\right)^k}{W_{mk+l}^3}\right)^{-1}-\left(-1\right)^n\left(W_{mn+l}^3+W_{m\left(n-1\right)+l}^3+I+J\right)\right)\\
		%=&\lim_{n\rightarrow \infty} \left(\left(-1\right)^nc_2^3\beta^{3\left(mn+l\right)}+\left(-1\right)^nc_2^3\beta^{3\left(m\left(n-1\right)+l\right)}+\frac{36\left(-1\right)^nc_2^3\beta^{3\left(mn+l\right)}\alpha^m\left(\alpha^{3m}-1\right)^3}{\left(\alpha^{4m}-\beta^{m}\right)\left(\alpha^{5m}-\beta^{2m}\right)}\right.\\
		%&\left.+\frac{36\left(-1\right)^nc_2^4\beta^{4\left(mn+l\right)}\alpha^m\left(\alpha^{3m}-1\right)^3}{c_1\alpha^{mn+l}\left(\alpha^{5m}-\beta^{2m}\right)^2} +O\left(\beta^{3\left(mn+l\right)}\right)\right)\\
		%=&\lim_{n\rightarrow \infty}\left(\left(-1\right)^nc_2^3\beta^{3\left(m\left(n-1\right)+l\right)}\left(\beta^{3m}+1+\frac{36\beta^{3m}\alpha^m\left(\alpha^{3m}-1\right)^3}{\left(\alpha^{4m}-\beta^{m}\right)\left(\alpha^{5m}-\beta^{2m}\right)}\right)\right.\\
		%&\left.+\left(\frac{\beta^4}{\alpha}\right)^{mn+l}\cdot \frac{36\left(-1\right)^nc_2^4\alpha^m\left(\alpha^{3m}-1\right)^3}{c_1\left(\alpha^{5m}-\beta^{2m}\right)^2}+O\left(\beta^{3\left(mn+l\right)}\right) \right)\\
		%=&0.
		%\end{aligned}
		%$$
		Therefore, by Definition \ref{definition 2.1}, we can get
		$$\left(\sum_{k=n}^{\infty}\frac{\left(-1\right)^n}{W_{mk+l}^3}\right)^{-1}\sim \left(-1\right)^n\left(W_{mn+l}^3+W_{m\left(n-1\right)+l}^3+I+J\right),$$
		where
		$$I=3c_1^2c_2\left(-\alpha q\right)^{mn+l}\left(1+\left(-\alpha q\right)^{-m}-\frac{\left(\alpha^{3m}+1\right)^2}{\alpha^{2m}\left(\alpha^{4m}+\beta^m\right)}\right)$$
		and 
		$$J=3c_1c_2^2\left(-\beta q\right)^{mn+l}\left(-1-\left(-\beta q\right)^{-m}-\frac{2\left(\alpha^{3m}+1\right)^2}{\alpha^m\left(\alpha^{5m}+\beta^m\right)}+\frac{3\left(\alpha^{3m}+1\right)^3}{\alpha^m\left(\alpha^{4m}+\beta^m\right)^2}\right),$$
		
		This completes the proof of Corollary \ref{corollary 4.3}.
		
		$\hfill\Box$
		
		\begin{corollary}\label{corollary 4.4}
			For any positive integers $m,n,p$ and integers $q,l$, where $l\geq 1-m$  and $p^2+2q-2<p\sqrt{p^2+4q}$, we have
			$$\left(\sum_{k=n}^{\infty}\frac{\left(-1\right)^k}{W_{mk+l}^4}\right)^{-1}\sim \left(-1\right)^n\left(W_{mn+l}^4+W_{m\left(n-1\right)+l}^4+L+M+N\right),$$
			where
			$$L=4c_1^3c_2\left(-\alpha^2 q\right)^{mn+l}\left(1+\left(-\alpha^2 q\right)^{-m}-\frac{\left(\alpha^{4m}+1\right)^2}{\alpha^{3m}\left(\alpha^{5m}+\beta^m\right)}\right),$$
			$$M=c_1^2c_2^2q^{2\left(mn+l\right)}\left(-6-6q^{-2m}-\frac{10\left(\alpha^{4m}+1\right)^2}{\alpha^{2m}\left(\alpha^{6m}+\beta^{2m}\right)}+\frac{16\left(\alpha^{4m}+1\right)^3}{\alpha^{2m}\left(\alpha^{5m}+\beta^m\right)^2}\right)$$
			and
			{\footnotesize$$
				N=4c_1c_2^3\left(-\beta^2 q\right)^{mn+l}\left(1+\left(-\beta^2 q\right)^{-m}-\frac{5\left(\alpha^{4m}+1\right)^2}{\alpha^m \left(\alpha^{7m}+\beta^{3m}\right)}-\frac{16\left(\alpha^{4m}+1\right)^4}{\alpha^m \left(\alpha^{5m}+\beta^m\right)^3}+\frac{20\left(\alpha^{4m}+1\right)^3}{c_1^2\alpha^m\left(\alpha^{5m}+\beta^m\right)\left(\alpha^{6m}+\beta^{2m}\right)}\right).
				$$}
			%$$\begin{aligned}
				%N=&4c_1c_2^3\left(-\beta^2 q\right)^{mn+l}\left(1+\left(-\beta^2 q\right)^{-m}-\frac{5\left(\alpha^{4m}+1\right)^2}{\alpha^m \left(\alpha^{7m}+\beta^{3m}\right)}-\frac{16\left(\alpha^{4m}+1\right)^4}{\alpha^m \left(\alpha^{5m}+\beta^m\right)^3}\right.\\
				%&\left.+\frac{20\left(\alpha^{4m}+1\right)^3}{c_1^2\alpha^m\left(\alpha^{5m}+\beta^m\right)\left(\alpha^{6m}+\beta^{2m}\right)}\right).
				%\end{aligned}$$
			\end{corollary}
			{\bf Proof}.\ By Theorem \ref{theorem 4.1}, we have
			$$\begin{aligned}
				&\left(\sum_{k=n}^{\infty}\frac{\left(-1\right)^k}{W_{mk+l}^4}\right)^{-1}\\
				=&\left(-1\right)^n\left(W_{mn+l}^4+W_{m\left(n-1\right)+l}^4\right)\\
				&-\left(-1\right)^n\sum_{i=1}^{4}\binom{4}{i}c_1^{4-i}\left(-c_2\right)^{i}\left(\left(\alpha^{mn+l}\right)^{4-i}\left(\beta^{mn+l}\right)^i+\left(\alpha^{m\left(n-1\right)+l}\right)^{4-i}\left(\beta^{m\left(n-1\right)+l}\right)^i\right)\\
				&+\frac{\left(-1\right)^{n}c_1^4\alpha^{4\left(mn+l\right)}\left(\alpha^{4m}+1\right)}{\alpha^{4m}}\sum_{j=1}^{3}\left(-1\right)^j\left(\sum_{i=1}^{3}\binom{3+i}{3}\frac{c_2^{i}\beta^{i\left(mn+l\right)}\alpha^{mi}\left(\alpha^{4m}+1\right)}{c_1^{i}\alpha^{i\left(mn+l\right)}\left(\alpha^{m\left(4+i\right)}+\beta^{mi}\right)}\right)^j+O\left(\beta^{4\left(mn+l\right)}\right)\\
				=&\left(-1\right)^n\left(W_{mn+l}^4+W_{m\left(n-1\right)+l}^4+L+M+N\right)+\left(-1\right)^nR,
			\end{aligned}$$
			where
			$$L=4c_1^3c_2\left(-\alpha^2 q\right)^{mn+l}\left(1+\left(-\alpha^2 q\right)^{-m}-\frac{\left(\alpha^{4m}+1\right)^2}{\alpha^{3m}\left(\alpha^{5m}+\beta^m\right)}\right),$$
			$$M=c_1^2c_2^2q^{2\left(mn+l\right)}\left(-6-6q^{-2m}-\frac{10\left(\alpha^{4m}+1\right)^2}{\alpha^{2m}\left(\alpha^{6m}+\beta^{2m}\right)}+\frac{16\left(\alpha^{4m}+1\right)^3}{\alpha^{2m}\left(\alpha^{5m}+\beta^m\right)^2}\right),$$
			{\footnotesize$$
				N=4c_1c_2^3\left(-\beta^2 q\right)^{mn+l}\left(1+\left(-\beta^2 q\right)^{-m}-\frac{5\left(\alpha^{4m}+1\right)^2}{\alpha^m \left(\alpha^{7m}+\beta^{3m}\right)}-\frac{16\left(\alpha^{4m}+1\right)^4}{\alpha^m \left(\alpha^{5m}+\beta^m\right)^3}+\frac{20\left(\alpha^{4m}+1\right)^3}{c_1^2\alpha^m\left(\alpha^{5m}+\beta^m\right)\left(\alpha^{6m}+\beta^{2m}\right)}\right)
				$$}
			and
			$$\begin{aligned}
				R=&-c_2^4\beta^{4\left(mn+l\right)}
				-c_2^4\beta^{4\left(m\left(n-1\right)+l\right)}+\frac{100c_2^4\beta^{4\left(mn+l\right)}\left(\alpha^{4m}+1\right)^3}{\left(\alpha^{6m}+\beta^{2m}\right)^2}+\frac{400c_2^6\beta^{6\left(mn+l\right)}\alpha^{2m}\left(\alpha^{4m}+1\right)^3}{c_1^2\alpha^{2\left(mn+l\right)}\left(\alpha^{7m}+\beta^{3m}\right)^2}\\
				&+\frac{160c_2^4\beta^{4\left(mn+l\right)}\left(\alpha^{4m}+1\right)^3}{\left(\alpha^{m}+\beta^{m}\right)\left(\alpha^{7m}+\beta^{3m}\right)}+\frac{400c_2^5\beta^{5\left(mn+l\right)}\alpha^{m}\left(\alpha^{4m}+1\right)^3}{c_1\alpha^{mn+l}\left(\alpha^{6m}+\beta^{2m}\right)\left(\alpha^{7m}+\beta^{3m}\right)}\\
				&-\frac{1000c_2^6\beta^{6\left(mn+l\right)}\alpha^{2m}\left(\alpha^{4m}+1\right)^4}{c_1^2\alpha^{2\left(mn+l\right)}\left(\alpha^{6m}+\beta^{2m}\right)^3}-\frac{8000c_2^9\beta^{9\left(mn+l\right)}\alpha^{5m}\left(\alpha^{4m}+1\right)^4}{c_1^5\alpha^{5\left(mn+l\right)}\left(\alpha^{7m}+\beta^{3m}\right)^3}\\
				&-\frac{480c_2^4\beta^{4\left(mn+l\right)}\left(\alpha^{4m}+1\right)^4}{\left(\alpha^{5m}+\beta^{m}\right)^2\left(\alpha^{6m}+\beta^{2m}\right)}-\frac{960c_2^5\beta^{5\left(mn+l\right)}\alpha^{m}\left(\alpha^{4m}+1\right)^4}{c_1\alpha^{mn+l}\left(\alpha^{5m}+\beta^{m}\right)^2\left(\alpha^{7m}+\beta^{3m}\right)}\\
				&-\frac{1200c_2^5\beta^{5\left(mn+l\right)}\alpha^{m}\left(\alpha^{4m}+1\right)^4}{c_1\alpha^{mn+l}\left(\alpha^{5m}+\beta^{m}\right)\left(\alpha^{6m}+\beta^{2m}\right)^2}-\frac{6000c_2^7\beta^{7\left(mn+l\right)}\alpha^{3m}\left(\alpha^{4m}+1\right)^4}{c_1^3\alpha^{3\left(mn+l\right)}\left(\alpha^{6m}+\beta^{2m}\right)^2\left(\alpha^{7m}+\beta^{3m}\right)}\\
				&-\frac{4800c_2^7\beta^{7\left(mn+l\right)}\alpha^{3m}\left(\alpha^{4m}+1\right)^4}{c_1^3\alpha^{3\left(mn+l\right)}\left(\alpha^{5m}+\beta^{m}\right)\left(\alpha^{7m}+\beta^{3m}\right)^2}-\frac{12000c_2^8\beta^{8\left(mn+l\right)}\alpha^{4m}\left(\alpha^{4m}+1\right)^4}{c_1^4\alpha^{4\left(mn+l\right)}\left(\alpha^{6m}+\beta^{2m}\right)\left(\alpha^{7m}+\beta^{3m}\right)^2}\\
				&-\frac{4800c_2^6\beta^{6\left(mn+l\right)}\alpha^{2m}\left(\alpha^{4m}+1\right)^4}{c_1^2\alpha^{2\left(mn+l\right)}\left(\alpha^{5m}+\beta^{m}\right)\left(\alpha^{6m}+\beta^{2m}\right)\left(\alpha^{7m}+\beta^{3m}\right)}.
			\end{aligned}$$
			Note that $p>0$, it's easy to prove that $|\alpha|>|\beta|$, furthermore, by  $p^2+2q-2<p\sqrt{p^2+4q}$, it's easy to prove that $|\beta|<1$. Therefore, for any positive integers $m$ and $n$ with $n<m$, we have $\left|\alpha^n\right|>\left|\beta^n\right|>\left|\beta^m\right|$, thus
			$$
			\begin{aligned}
				&\lim_{n\rightarrow \infty}
				\left(\left(\sum_{k=n}^{\infty}\frac{\left(-1\right)^k}{W_{mk+l}^4}\right)^{-1}-\left(-1\right)^n\left(W_{mn+l}^4+W_{m\left(n-1\right)+l}^4+L+M+N\right)\right)\\
				=&\lim_{n\rightarrow \infty}
				\left(\left(-1\right)^nR+O\left(\beta^{4\left(mn+l\right)}\right)\right)\\
				=&0.
			\end{aligned}
			$$
			Therefore, by Definition \ref{definition 2.1}, we have
			$$\left(\sum_{k=n}^{\infty}\frac{\left(-1\right)^k}{W_{mk+l}^4}\right)^{-1}\sim \left(-1\right)^n\left(W_{mn+l}^4+W_{m\left(n-1\right)+l}^4+L+M+N\right),$$
			where
			$$L=4c_1^3c_2\left(-\alpha^2 q\right)^{mn+l}\left(1+\left(-\alpha^2 q\right)^{-m}-\frac{\left(\alpha^{4m}+1\right)^2}{\alpha^{3m}\left(\alpha^{5m}+\beta^m\right)}\right),$$
			$$M=c_1^2c_2^2q^{2\left(mn+l\right)}\left(-6-6q^{-2m}-\frac{10\left(\alpha^{4m}+1\right)^2}{\alpha^{2m}\left(\alpha^{6m}+\beta^{2m}\right)}+\frac{16\left(\alpha^{4m}+1\right)^3}{\alpha^{2m}\left(\alpha^{5m}+\beta^m\right)^2}\right)$$
			and
			{\footnotesize$$
				N=4c_1c_2^3\left(-\beta^2 q\right)^{mn+l}\left(1+\left(-\beta^2 q\right)^{-m}-\frac{5\left(\alpha^{4m}+1\right)^2}{\alpha^m \left(\alpha^{7m}+\beta^{3m}\right)}-\frac{16\left(\alpha^{4m}+1\right)^4}{\alpha^m \left(\alpha^{5m}+\beta^m\right)^3}+\frac{20\left(\alpha^{4m}+1\right)^3}{c_1^2\alpha^m\left(\alpha^{5m}+\beta^m\right)\left(\alpha^{6m}+\beta^{2m}\right)}\right).
				$$}
			
			This completes the proof of Corollary \ref{corollary 4.4}.
			
			$\hfill\Box$
			
			\section{Conclusions}
			%In this paper, for the generalized Fibonacci sequence $\left\{W_n\left(a,b,p,q\right)\right\}$, we consider the asymptotic estimation for two classes of generalized Fibonacci sub-sequences and obtain the following main results, which generalize  the asymptotic estimation results of Yuan et al. \cite{A22} in 2025.
			In this paper, we consider the asymptotic estimation for  generalized Fibonacci sub-sequences $\left\{W_n\left(a,b,p,q\right)\right\}$ and obtain the following two main results.
			\begin{itemize}
				\item The asymptotic estimate of $\left(\sum\limits_{k=n}^{\infty}\frac{1}{W_{mk+l}^d}\right)^{-1}$(Theorem \ref{theorem 3.1}). In particular, for some special $a,b,p,q,m,l$ and $d$, we have the following table.
				\begin{center}\renewcommand{\arraystretch}{1.5}
					\begin{tabular}{|c|c|c|c|c|c|} % 每列宽度1.2cm，内容居中
						\hline
						 $\left(a,b\right) \in \mathbb{Z}^2$ & $l \in \mathbb{N}$ & $\left(p,q\right) \in \mathbb{Z}^2$  &  $m\in \mathbb{Z}^+$& $d\in \mathbb{Z}^+$ & Theorem \ref{theorem 3.1} \\
						\hline
						 \multirow{10}{*}{$\left(0,1\right)$} &  \multirow{10}{*}{0} &\multirow{6}{*}{$\left(1,1\right)$} & 1 & 2  & \cite{A11}Theorem 2.1 \\ 
						\cline{4-6}
						  &   & &3 & 2  &  \cite{A11}Theorem 5.1 \\ 
						\cline{4-6}
						 &   & &1 & 1  &  \cite{A12}Theorem 1 \\ 
						\cline{4-6}
						  &   & &$\mathbb{Z}^+$ & 1  & \cite{A12}Theorem 2 \\ 
						\cline{4-6}
						  &   & &$\mathbb{Z}^+$ & 2  &  \cite{A13}Theorem 1 \\ 
						\cline{4-6}
						  & &  & 1& 4  &  \cite{A10}Corollary 1 \\ 
						\cline{3-6}
						  &   &\multirow{4}{*}{$\left(\mathbb{Z}^+,\pm1\right)$}  & \multirow{4}{*}{$\mathbb{Z}^+$} & 1   &\cite{A14}Theorem 2.1 \\ 
						\cline{5-6}
						& & & & 2 &   \cite{A14}Theorem 2.2 \\ 
						\cline{5-6}
						&& & & 3 &  \cite{A14}Theorem 2.3 \\ \cline{5-6}
						 & &&&4&  \cite{A14}Theorem 2.4 \\  
						\hline
					\end{tabular}
				\end{center}
				\item The asymptotic estimate of $\left(\sum\limits_{k=n}^{\infty}\frac{\left(-1\right)^k}{W_{mk+l}^d}\right)^{-1}$ (Theorem \ref{theorem 4.1}), which is the first result for the asymptotic estimate of $\left(\sum\limits_{k=n}^{\infty}\frac{\left(-1\right)^k}{W_{mk+l}^d}\right)^{-1}$ so far.
			\end{itemize}

			%In this paper, for the generalized Fibonacci sequence $\left\{W_n\left(a,b,p,q\right)\right\}$, we first establish the corresponding connections of $\left(\sum\limits_{k=n}^{\infty}\frac{1}{W_{mk+l}^d}\right)^{-1}$ or $\left(\sum\limits_{k=n}^{\infty}\frac{\left(-1\right)^k}{W_{mk+l}^d}\right)^{-1}$ and their asymptotic estimations, respectively, and then obtain the exact asymptotic estimation values of $\left(\sum\limits_{k=n}^{\infty}\frac{1}{W_{mk+l}^d}\right)^{-1}$ or $\left(\sum\limits_{k=n}^{\infty}\frac{\left(-1\right)^k}{W_{mk+l}^d}\right)^{-1}$ for $d=1,2,3,4$, respectively. In particular, we generalize  the asymptotic estimation results of Yuan et al. in \cite{A14} as the following table.
			
			%In this paper, we consider the asymptotic estimation for two classes of generalized Fibonacci sub-sequences and obtain the following main results, which generalize  the asymptotic estimation results of Yuan et al. \cite{A14} in 2025. 
			%\begin{itemize}\item The asymptotic estimate of $\left(\sum\limits_{k=n}^{\infty}\frac{1}{W_{mk+l}^d}\right)^{-1}$(Theorem \ref{theorem 3.1}). In particular, we have the following table.
			%\item The asymptotic estimate of $\left(\sum\limits_{k=n}^{\infty}\frac{\left(-1\right)^k}{W_{mk+l}^d}\right)^{-1}$(Theorem \ref{theorem 4.1}).
			%\end{itemize}
			
			\section*{Acknowledgements}
			This paper is supported by National Natural Science Foundation of China(Grant No. 12471494) and Natural Science Foundation of Sichuan Province (2024NSFSC2051). The corresponding author is Professor Qunying Liao.
			
\end{document}